\documentclass[12pt,reqno]{amsart}
\usepackage{graphicx}
\usepackage{epstopdf}
\usepackage{amsmath}
\usepackage{amsfonts}
\usepackage{amssymb}
\usepackage[parfill]{parskip}
\usepackage{setspace}
\usepackage[small,bf,hang]{caption}
\usepackage{subcaption}
\usepackage{float}
\usepackage{varioref}
\usepackage{algorithmic}
\usepackage[numbers,sort]{natbib}
\usepackage[subnum]{cases}

\usepackage[
  hmarginratio={1:1},     
  vmarginratio={1:1},     
  textwidth=6in,        
  heightrounded,          
]{geometry}
\usepackage{etoolbox}
\patchcmd{\section}{\scshape}{\bfseries}{}{}
\patchcmd{\subsection}{\normalfont}{\itshape}{}{}
\makeatletter
\renewcommand{\@secnumfont}{\bfseries}
\makeatother

\labelformat{equation}{Eq.~#1}
\labelformat{figure}{Figure~#1}
\labelformat{table}{Table~#1}

\newcommand{\half}{\frac{1}{2}}
\newcommand{\pd}[2]{\frac{\partial #1}{\partial #2}}

\newcommand{\ang}[1]{\langle #1 \rangle}

\newcommand{\mb}[1]{\textbf{\textit{#1}}}

\DeclareMathOperator{\sgn}{sgn}

\newcommand{\tab}{\hskip15pt}
\newcommand{\sfwidth}{0.475\textwidth}
\newcommand{\inputImgEps}[1] {\includegraphics[width=0.9\textwidth]{#1}}

\newcommand{\phalf}{\mb{i} + \half \mb{e}^d}
\newcommand{\mhalf}{\mb{i} - \half \mb{e}^d}
\newcommand{\pmhalf}{\mb{i} \pm \half \mb{e}^d}
\newcommand{\splus}{\mb{i} + \mb{e}^d}
\newcommand{\sminus}{\mb{i} - \mb{e}^d}

\title[Single Stage Flux-Corrected Transport Algorithm]{A Single Stage Flux-Corrected Transport Algorithm for High-Order Finite-Volume Methods}

\author{Christopher Chaplin} 
\author{Phillip Colella}

\begin{document}

\begin{abstract}
We present a new limiter method for solving the advection equation using a high-order, finite-volume discretization.
The limiter is based on the flux-corrected transport algorithm.
We modify the classical algorithm by introducing a new computation for solution bounds at smooth extrema, as well as improving the pre-constraint on the high-order fluxes.
We compute the high-order fluxes via a method of lines approach with fourth order Runge-Kutta as the time integrator.
For computing low-order fluxes, we select the corner transport upwind method due to its improved stability over donor-cell upwind.
Several spatial differencing schemes are investigated for the high-order flux computation, including centered difference and upwind schemes.
We show that the upwind schemes perform well on account of the dissipation of high wavenumber components.
The new limiter method retains high-order accuracy for smooth solutions and accurately captures fronts in discontinuous solutions.
Further, we need only apply the limiter once per complete time step.
\end{abstract}

\maketitle

\section{Introduction}

We wish to solve hyperbolic conservation laws of the form

\begin{equation} \label{consv}
\pd{U}{t} + \nabla \cdot ( \vec{\mb{F}}(U) ) = 0
\end{equation}

where $U$ represents a vector of conserved values and $\vec{\mb{F}}(U), [\vec{\mb{F}} = (\mb{F}^1 \dots \mb \mb{F}^D)]$ are corresponding fluxes. 
The discrete solution of these equations at a given time $t^{n+1}$ and spatial location $\mb{i}$ is given by:

\begin{equation} \label{dcl}
\ang{U}^{n+1}_{\mb{i}} = \ang{U}^{n}_{\mb{i}} - \frac{\Delta t}{\Delta x} \sum_d \left[ (\mb{F}^d)^{n+\half}_{\phalf} - (\mb{F}^d)^{n+\half}_{\mhalf} \right]
\end{equation}

where $\ang{U}^n$ and $(\mb{F}^d)^{n+\half}$ approximate the average of $U$ over a rectangular Cartesian control volumes
indexed by $\mb{i}$ at time step $n$ and the average of the fluxes in time over the respective faces, respectively. 
 


Methods for accurately computing these fluxes $(\mb{F}^d)^{n+\half}$ to obtain high-order accuracy for smooth solutions are well-understood. 
However, in the presence of discontinuities or underresolved gradients,
such methods generate oscillatory errors, and limiters are typically introduced to selectively introduce dissipation that will damp these 
errors. For high-order methods based on the method of lines, the typical approach to limiters has been to apply a limiter independently at every flux evaluation in a multi-stage high-order time integration method, e.g. Runge-Kutta methods. However, such an approach is problematic, particularly for limiters that are
intended to preserve high-order accuracy at smooth extrema. Such limiters typically make use of second- or higher-order spatial derivative information at extrema to determine whether the solution is smooth enough to not require limiting. This leads to the use of a larger stencil (for the limiter) at every stage of the time integration than that of the simpler high-order spatial discretization of the fluxes without limiting, thus increasing communication costs in a domain-decomposed parallel algorithm. In addition, only the linear combination of Runge-Kutta stage computations is guaranteed to produce a high-order accurate solution; the individual stage fluxes and solutions are typically lower-order accurate in time. In adaptive mesh refinement methods with refinement in time, this feature complicates the design of the time interpolation step for computing ghost cell values at refinement boundaries \cite{Colella2011}. Finally, many of the stage-by-stage limiting schemes are geometric in nature, such that the effective fully limited method corresponds to donor-cell spatial differencing for the fluxes \cite{harten1987,liu1994}. This typically has a more restrictive CFL time step stability condition than the high-order methods. 

The starting point for addressing these problems is to use a version of the classic Flux-Corrected Transport (FCT) algorithms \cite{Boris1973,Zalesak1979,FCT_Book}.
These methods introduce dissipation through a nonlinear hybridization of a high-order flux with a dissipative, low-order flux. 
The FCT algorithm is particularly advantageous because it has a straightforward multi-dimensional expression \cite{Zalesak1979}
and has been used extensively with finite-volume discretizations. We will look at a family of high-order methods based on centered- and one-point-upwinded linear finite volume interpolation, following the ideas in \cite{Colella2011}, combined with the method of lines with fourth-order Runge-Kutta time discretization as the method for time discretization. All four intermediate stages of this high-order method will be computed without any limiting applied, in order to produce a high-order flux to go into the FCT hybridization. The low-order method will be the corner-transport upwind (CTU) method 
\cite{Colella1990,Saltzman1994}. The method for computing the hybridization coefficient will be based on the approach used in \cite{Colella2008,Colella2011} for interpolation-based limiting, modified to fit into the formalism in \cite{Zalesak1979}. This combination of methods addresses the issues alluded to above. The CFL stability condition for the CTU method is independent of the direction of propagation of the waves relative to the coordinate directions and of the dimensionality of the problem, and is typically much less restrictive than the CFL condition for donor-cell / method-of-lines. The limiter is applied only once -- at the end of the time step -- thus minimizing the impact on the communications cost. 

For this study, we restrict our attention to the equation of scalar advection equation. This allows us to explore design decisions in a simple setting, but one that is still relatively unforgiving, and important for real applications (transport of scalars in the atmosphere, Vlasov equations in phase space). The extension to systems will be discussed in future work. One question we will not address here is that of positivity preservation. The reason is that, in more than one dimension for any linear first-order method that is not donor cell, it is possible to construct a discretely divergence-free advection velocity and non-negative initial data such that the solution becomes negative after one time step. Our preferred approach is to post-process the solution at the end of each time step by redistributing negative increments to nearby cells in a way that would lead to an overall non-negative solution, while still conserving \cite{Hilditch1995,Wang2011}.

\subsection{Advection Equation}
We will consider the linear advection equation in the following form.

\begin{align} \label{advec}
\pd{q}{t} + \nabla \cdot (q\vec{\mb{u}}) &= 0 \\
\nabla \cdot \vec{\mb{u}} &= 0
\end{align}

on a $D$-dimensional square domain $\Omega = [0,1]^D$.
In this case $\vec{\mb{u}}$ is an advective velocity and $q$ is a scalar field.
The partial differential equation above can also be written as:

\begin{equation} \label{advec_solution}
\frac{dq}{dt} = 0 \tab \frac{d\vec{\mb{x}}}{dt} = \vec{\mb{u}}
\end{equation}

Provided that an initial condition is specified ($q_0 = q(\vec{\mb{x}}(t_0),t_0)$) this system of ordinary differential equations yields a unique solution for any $q(\vec{\mb{x}}(t),t)$ and $\vec{\mb{x}}(t)$.
The solution arrived at by integrating the equations is that $q$ is constant along characteristic curves defined by $\vec{\mb{x}}(t)$.
Even though there is a simple solution to this equation, the analysis is still quite useful since there is no diffusion or entropy condition built-in to the equation:
any numerical errors introduced are propagated through the domain.

\subsection{Finite-Volume Discretization}
Our solution approach is to use a finite-volume method to discretize the physical domain into a union of control volumes $V_{\mb{i}}$ (\ref{cv})

\begin{equation} \label{cv}
V_{\mb{i}} = [\mb{i}h , (\mb{i}+\mb{e})h] \ , \ \mb{i} \in {\mathbb{Z}}^{D} \ , \ \mb{e} = (1, 1, ..., 1)
\end{equation}

where $h$ is the grid spacing ($\Delta x$) and $\mb{i}$ is a $D$-dimensional index denoting location.
Values of the conserved scalar quantity $q$ are stored as cell-averages $\ang{q}$ over each cell $V_{\mb{i}}$ (\ref{avg}).
The fluxes $\mb{F}^d = q\mb{u}^d$ are stored as averages $\ang{\mb{F}^d}_{\pmhalf}$ over the surface faces $A_{d}^{\pm}$ of each cell (\ref{favg}).

\begin{equation} \label{avg}
\ang{q}_{\mb{i}}(t) = \frac{1}{h^D} \int_{V_{\mb{i}}} q(\mb{x},t) d\mb{x}
\end{equation}
\begin{equation} \label{favg}
\ang{\mb{F}^d}_{\pmhalf} (t) = \frac{1}{h^{D-1}} \int_{A_{d}^{\pm}} \mb{F}^d(\mb{x},t) \ d\mb{x}
\end{equation}

Applying the finite-volume discretization (\ref{cv}) to \ref{advec} yields a semi-discrete system of ordinary differential equations (ODE) in time:

\begin{align}
\label{dpde}
\frac{d\ang{q}_{\mb{i}}}{dt} &= - \frac{1}{h^D} \int_{V_{\mb{i}}} ( \nabla \cdot (\vec{\mb{F}}) ) d\mb{x} \\
\label{ddiv}
			                 &= -\frac{1}{h} \sum_d [ \ang{\mb{F}^d}_{\phalf} - \ang{\mb{F}^d}_{\mhalf} ]
\end{align}

where \ref{ddiv} is the result of applying the divergence theorem to \ref{dpde}.
The integration of the above system with respect to time from $t^n$ to $t^{n+1}$ produces the solution given below:

\begin{align} \label{fv_dcl}
\ang{q}^{n+1}_{\mb{i}} &= \ang{q}^{n}_{\mb{i}} - \frac{\Delta t}{h} \sum_d \left[ \ang{\mb{F}^d}^{n+\half}_{\phalf} - \ang{\mb{F}^d}^{n+\half}_{\mhalf} \right] \\
\label{fv_avg_flux}
\ang{\mb{F}^d}^{n+\frac{1}{2}}_{\pmhalf} &= \frac{1}{\Delta t} \int_{t^n}^{t^n+\Delta t} \ang{\mb{F}^d}_{\pmhalf} (t) dt
\end{align}

The resulting challenge is to accurately compute $\ang{\mb{F}^d}^{n+\half}_{\pmhalf}$.
It is important to note that no approximations have been made at this point: \ref{fv_dcl} - \ref{fv_avg_flux} are exact relationships.
However, to obtain a full discrete approximation, we need quadrature rules for the surface fluxes in \ref{ddiv} and for the time-averaged fluxes in \ref{fv_avg_flux}.
The quadrature rules for computing these fluxes are defined following ideas from \cite{Colella2011}.
In that work, the high-order quadratures were computed using a method of lines approach.
The surface fluxes were computed using a high-order centered difference method and the temporal integration was computed using the classic fourth-order Runge-Kutta (RK4) method.
We retained the use of RK4 in this study and investigated several high-order methods for computing the surface fluxes.

\subsection{Hybridization}
Returning to the flux description in \ref{fv_avg_flux}, we may now define the hybridization

\begin{equation} \label{hybrid}
\ang{\mb{F}^d}^{n+\half}_{\phalf} = (\eta_{\phalf}) \ang{\mb{F}^d_H}_{\phalf} + (1-\eta_{\phalf}) \ang{\mb{F}^d_L}_{\phalf}
\end{equation}

where the subscripts $H$ and $L$ refer to the high-order and low-order fluxes,
and $\eta_{\phalf}$ is the hybridization coefficient.

In the following sections of the paper we will describe the design choices and procedures for computing
the high-order flux, the low-order flux, and the hybridization coefficient.

\section{High-Order Flux Computation}
We compute the high-order fluxes using the method of lines.
Two schemes must be chosen: a scheme for integrating the solution in time and a scheme for computing the spatial derivatives.
High-order accuracy requires that both schemes be high-order accurate.

\subsection{High-Order Temporal Integration Scheme}
We use the RK4 scheme to advance the solution forward in time.
Returning to the system of ODE (\ref{dpde}):

\begin{align} \label{dcl2}
\frac{d\ang{q}}{dt} &= - D \cdot \ang{\vec{\mb{F}}} (t) \\
                             &= -\frac{1}{h} \sum_d \ang{\mb{F}^d}_{\phalf} - \ang{\mb{F}^d}_{\mhalf} \nonumber
\end{align}

we want to integrate the system from $t^n$ to $t^{n+1}$.
RK4 is a fourth-order integration scheme that consists of computing a linear combination of stage update variables $k_s$.
The updates are defined below:

\begin{alignat}{2}
\ang{q}^{0} &= \ang{q} (t^n) \tab &k_1 = - D \cdot \ang{\vec{\mb{F}} ( \ang{q}^0)  } \Delta t  \\ 
\ang{q}^1 &=  \ang{q}^0 + \frac{k_1}{2} \tab  &k_2 = - D \cdot \ang{\vec{\mb{F}} ( \ang{q}^1)} \Delta t \\ 
\ang{q}^2 &= \ang{q}^0 + \frac{k_2}{2} \tab &k_3 = - D \cdot \ang{\vec{\mb{F}} ( \ang{q}^2)} \Delta t \\
\ang{q}^3 &= \ang{q}^0 + k_3 \tab  &k_4 = - D \cdot \ang{\vec{\mb{F}} (\ang{q}^3)} \Delta t 
\end{alignat}

Each update variable $k_s$ requires computing stage fluxes $\ang{\mb{F}^d}^{s}_{\pmhalf} = \ang{q\mb{u}^d}^{s}_{\pmhalf}$.
The stage fluxes are functions of the stage values $\ang{q}^s_{\mb{i}}$ and $\ang{\mb{u}^d}^s_{\mb{i}}$ alone and the procedure for computing the fluxes will be described in the next section.

To perform the RK4 integration:

\begin{equation}
\ang{q} (t^n + \Delta t) = \ang{q} (t^n) + \frac{1}{6} (k_1 + 2k_2 + 2k_3 + k_4)
\end{equation}

Returning to the conservation notation, this RK4 integration can be described by 

\begin{align}
\ang{q}^{n+1}_{\mb{i}} &= \ang{q}^n_{\mb{i}} - \frac{\Delta t}{h} \sum_{d=1}^D [ \ang{\mb{F}^d_H}_{\phalf} - \ang{\mb{F}^d_H}_{\mhalf} ] \\
\label{highFlux}
\ang{\mb{F}^d_H}_{\pmhalf} &= \frac{1}{6} \left[\ang{\mb{F}^d}^{(0)}_{\pmhalf} + 2\ang{\mb{F}^d}^{(1)}_{\pmhalf} 
					+ 2\ang{\mb{F}^d}^{(2)}_{\pmhalf} + \ang{\mb{F}^d}^{(3)}_{\pmhalf}\right]
\end{align}

\subsection{High-Order Spatial Difference Schemes}
\begin{figure} [tb]
\centering
	\begin{subfigure}[b]{\sfwidth}
		\resizebox{\textwidth}{!}{\input{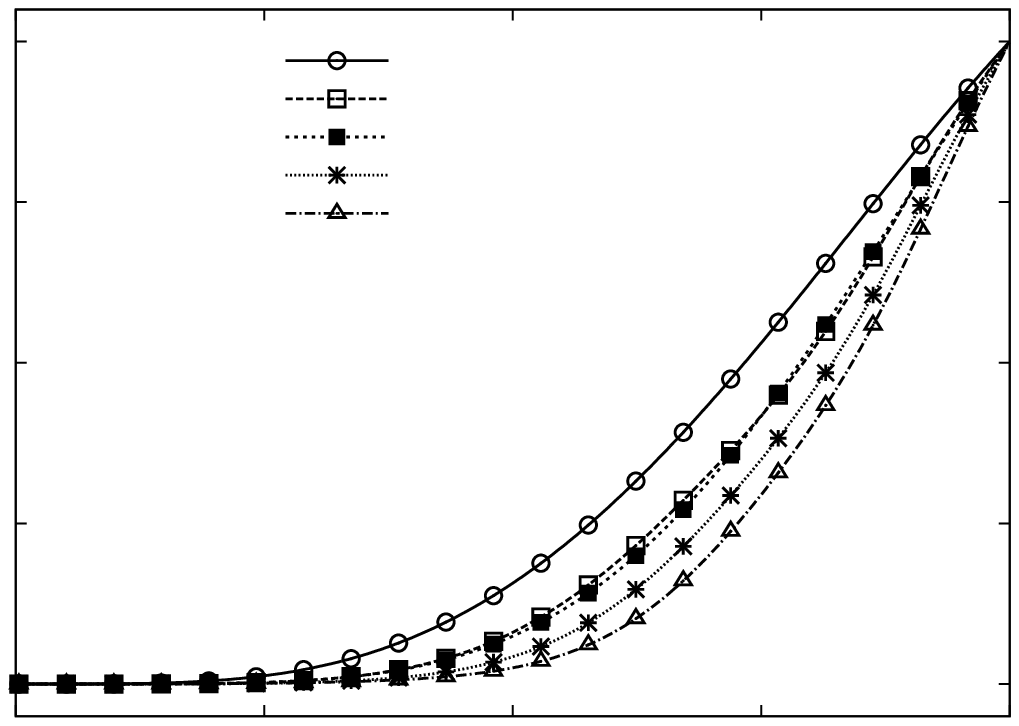}}
		\caption{Phase Error}
		\label{phase}
	\end{subfigure}
	\begin{subfigure}[b]{\sfwidth}
		\resizebox{\textwidth}{!}{\input{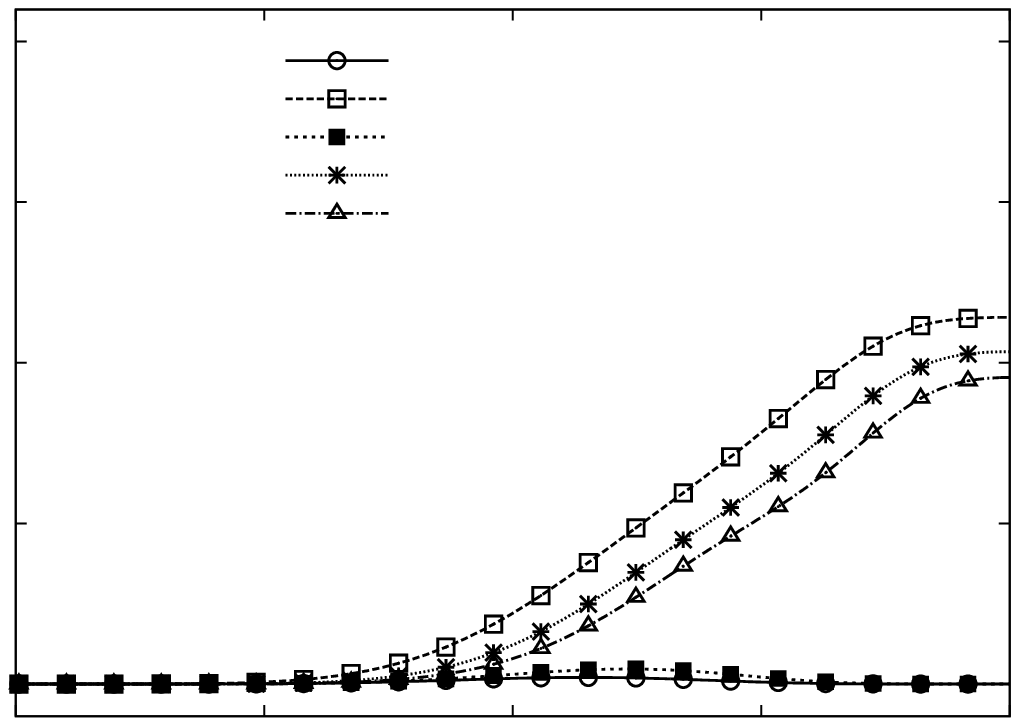}}
		\caption{Dissipation}
		\label{diss}
	\end{subfigure}
\captionsetup{width=\textwidth}
\caption[phase and diss]{Normalized phase error and dissipation for the high-order methods ($\sigma$ = 0.8)}
\label{phase and diss}
\end{figure}

We use high-order finite-difference methods to approximate the surface fluxes associated with the spatial derivatives.
The fluxes $\ang{q\mb{u}^d}_{\pmhalf}$ for the spatial derivatives are functions only of the cell-averaged $\ang{q}_{\mb{i}}$ and $\ang{\mb{u}^d}_{\mb{i}}$ 
at any time.
Several methods were explored in this study for computing $\ang{q}_{\pmhalf}$, including high-order centered difference schemes and upwind schemes. 
The advantage of the upwind methods is that they have greater diffusion especially in regimes where the phase error begins to rise (\ref{phase and diss}).
The upwind methods only require a small additional computation and the stability between similar order centered difference and upwind methods is almost identical (\ref{stab}).

The interpolation formulae corresponding to the spatial differencing schemes used are presented below.
For compactness, the following notation will be used:

\begin{equation}
\ang{q}^{n}_{\phalf} = \sum^{S}_{s=-S} a_s \ang{q}^n_{\mb{i} + s\mb{e}^d}
\end{equation}

where $S$ is the width of the stencil and $a_s$ are the coefficients. 
The odd ordered methods use the full range of coefficients, whereas the even ordered methods have no coefficient at $s = -S$.

\emph{Fourth Order Centered Difference} $(S = 2)$
\begin{equation}
\left\{ a_s : s = -S+1,\ldots,S\right\} = \frac{1}{12} \left\{ -1,7,7,-1 \right\}
\end{equation}
\emph{Fifth Order Upwind} $(S = 2)$
\begin{equation}
\left\{ a_s : s = -S,\ldots,S\right\} = \frac{1}{60} \left\{ 2,-13,47,27,-3 \right\}
\end{equation}
\emph{Sixth Order Centered Difference} $(S = 3)$
\begin{equation}
\left\{ a_s : s = -S+1,\ldots,S\right\} = \frac{1}{60} \left\{ 1,-8,37,37,-8,1 \right\}
\end{equation}
\emph{Seventh Order Upwind} $(S = 3)$
\begin{equation}
\left\{ a_s : s = -S,\ldots,S\right\} = \frac{1}{420} \left\{ -3,25,-101,319,214,-38,4 \right\}
\end{equation}
\emph{Ninth Order Upwind} $(S = 4)$
\begin{equation}
\left\{ a_s : s = -S,\ldots,S\right\} = \frac{1}{2520} \left\{ 4,-41,199,-641,1879,1375,-305,55,-5 \right\}
\end{equation}

\subsection{Product Rule}
To complete the flux computation, we must compute the average of the product of the scalar variable and the velocity ($\ang{q\mb{u}^d}_{\phalf}$).
The 2D product rules for second, fourth, and sixth-order accuracy are:

\begin{align}
\ang{q\mb{u}^d}_{\phalf} = &\ang{q}_{\phalf}\ang{\mb{u}^d}_{\phalf} + O(h^2) \\
\ang{q\mb{u}^d}_{\phalf} = &\ang{q}_{\phalf}\ang{\mb{u}^d}_{\phalf} + \frac{h^2}{12}\sum_{d' \ != \ d}{ \pd{q}{x^{d'}} \pd{\mb{u}^d}{x^{d'}} } + O(h^4) \\
\ang{q\mb{u}^d}_{\phalf} = &\ang{q}_{\phalf}\ang{\mb{u}^d}_{\phalf} + \frac{h^2}{12}\sum_{d' \ != \ d}\left({ \pd{q}{x_{d'}} \pd{\mb{u}^d}{x_{d'}} }\right) \\ \nonumber
				& + \frac{h^4}{1440}\sum_{d' \ != \ d}{\left(3\frac{\partial^3 q}{\partial x^3_{d'}}\pd{\mb{u}^d}{x_{d'}} +
							3\frac{\partial^3 \mb{u}^d}{\partial x^3_{d'}}\pd{q}{x_{d'}} + 
							2\frac{\partial^2 \mb{u}^d}{\partial x^2_{d'}}\frac{\partial^2 q}{\partial x^2_{d'}}\right)}
							+ O(h^6)
\end{align}

The possible sources of error in the product formulae above are computing the averages $\ang{q}_{\phalf}$, $\ang{\mb{u}^d}_{\phalf}$, and computing the partial derivative sums.
We have already discussed several methods and their accuracy for computing $\ang{q}_{\phalf}$.
The velocity fields are analytic for advection, so $\ang{\mb{u}^d}_{\phalf}$ introduces no error.
The derivative terms in the summations above were computed exclusively using centered difference approximations of appropriate accuracy.
For example, the derivatives in the fourth-order accurate product formula were computed using a second-order centered difference.
The derivatives in the sixth-order formula were computed to fourth-order (for the term multiplied by $h^2$) and to second-order (for the term multiplied by $h^4$).

\subsection{Stability}
We compute the stability for each high-order scheme to determine the allowable time step size following the procedure in \cite{Colella2011a}.
Stability for the method of lines requires the eigenvalues of the right hand side to lie within the stability region of the time integrator.
These eigenvalues are computed by diagonalizing the the semi-discrete system (\ref{dcl2}).
For advection the eigenvalues are defined as the product of the velocity and the spatial derivative operator (\ref{MOL}-\ref{eigen}).
The particular eigenvalues for each spatial differencing scheme will be presented later.

\begin{align} \label{MOL}
\frac{d\ang{q}}{dt} &= \lambda \ang{q} \\
\label{eigen}
\lambda \ang{q} &= -\vec{\mb{u}} \pd{}{\mb{x}} \ang{q}
\end{align}

These eigenvalues must lie within the stability region of the time integrator.
The stability region for RK4 is well known, and can be described by its characteristic polynomial

\begin{equation}\label{poly}
P(z) = 1 + z + \frac{z^2}{2} + \frac{z^3}{6} + \frac{z^4}{24}
\end{equation}

where $z = \Delta t \lambda$.
Stability for this problem requires that $|P(z)| \le 1$.
The resulting stability constraints for each spatial differencing scheme are presented in \ref{stab},
where $\sigma = |u|\Delta t/h$.

\begin{table}[tb]
\caption{Stability of Method for Varying Spatial Difference Operators} \label{stab}
\centering
    \begin{tabular}{ | p{3.00cm} | p{5.00cm} |}
    \hline
    \bf{Method} & \bf{Stability Constraint} \\ \hline \hline
    4th Center & $\sigma \lesssim $ 2.06/D  \\ \hline
    5th Upwind & $\sigma \lesssim $ 1.73/D \\ \hline
    6th Center & $\sigma \lesssim $ 1.78/D   \\ \hline
    7th Upwind & $\sigma \lesssim $ 1.69/D    \\ \hline
    9th Upwind & $\sigma \lesssim $ 1.60/D   \\ \hline
    \end{tabular}
\end{table}

Along with stability, the phase error and dissipation were computed (\ref{phase and diss}).
The dissipation was defined as $(1 - |g|)$ where:

\begin{align}
|g| &= \sqrt{\left(Re(g)^2 + Im(g)^2\right)} \\
Re(g) &= \left( 1 + x + \frac{x^2}{2} + \frac{x^3}{6} + \frac{x^4}{24} \right) - \frac{y^2}{2} \left(1 + x + \frac{x^2}{2} \right) + \frac{y^4}{24}   \\
Im(g) &= y \left( 1 + x + \frac{x^2}{2} + \frac{x^3}{6} \right) - \frac{y^3}{6}\left( 1+ x \right)
\end{align}

and $\lambda = x + i y$.
The normalized phase error, $|1-\alpha|$, is defined using the following computation

\begin{equation}
\alpha = \frac{\alpha (\beta)}{|u|} = -\frac{1}{\sigma \beta} \frac{Im(g)}{Re(g)}
\end{equation}

\subsubsection{Spatial Differencing Eigenvalues}
The eigenvalues for each of the different high-order spatial differencing schemes are presented below.
In each of the eigenvalue descriptions, $\beta_d$ may range from $-\pi$ to $\pi$, and is defined as $2\pi k_d h$ with $k_d = 0,\pm1 ,\pm2,...,\pm N/2$

\emph{Fourth Order Centered Difference}:
\begin{equation} \label{fourth}
\lambda_4 = \frac{i}{12h} \sum_{d=1}^D u^d \left[ 16\sin(\beta_d) - 2\sin(2\beta_d) \right]
\end{equation}

\emph{Fifth Order Upwind}:
\begin{align} \label{fifth}
\lambda_5 = \frac{1}{60h} \sum_{d=1}^D u^d [ &\left( -2\cos(3\beta_d) + 12\cos(2\beta_d) - 30\cos(\beta_d) + 20 \right) + \nonumber \\
						  &i\left( 2\sin(3\beta_d) - 18\sin(2\beta_d) + 90\sin(\beta_d) \right) ]
\end{align}

\emph{Sixth Order Centered Difference}:
\begin{equation} \label{sixth}
\lambda_6 = \frac{i}{60h} \sum_{d=1}^D u^d \left[ 2\sin(3\beta_d) - 18\sin(2\beta_d) + 90\sin(\beta_d) \right]
\end{equation}

\emph{Seventh Order Upwind}:
\begin{align} \label{seventh}
\lambda_7 = \frac{1}{420h} \sum_{d=1}^D u^d [&\left( 3\cos(4\beta_d) - 24\cos(3\beta_d) + 84\cos(2\beta_d) - 168\cos(\beta_d) + 105 \right) + \nonumber \\
						  &i\left( -3\sin(4\beta_d) + 32\sin(3\beta_d) - 168\sin(2\beta_d) + 672\sin(\beta_d) \right) ]
\end{align}

\emph{Ninth Order Upwind}:
\begin{align} \label{ninth}
\lambda_9 = \frac{1}{2520h} \sum_{d=1}^D u^d [&( -4\cos(5\beta_d)+40\cos(4\beta_d)-180\cos(3\beta_d)+ \nonumber \\
					     &480\cos(2\beta_d)-840\cos(\beta_d)+504) + \nonumber \\
					     &i( 4\sin(5\beta_d)-50\sin(4\beta_d)+300\sin(3\beta_d)- \nonumber \\
					     &1200\sin(2\beta_d)+4200\sin(\beta_d)) ]
\end{align}

\section{Low-Order Flux Computation}
The low order fluxes are computed using the CTU method \cite{Colella1990,Saltzman1994}.
CTU is a first-order time advancement scheme.
The method is desirable over the simpler donor-cell upwind method because its
stability is independent of dimensionality.
However, this increased stability comes with a price.
Instead of a single flux being defined by a single upwind value, the CTU flux is dependent upon
a set of upwinded values.
These values are determined by tracing the characteristic paths from the nodes that define the flux surface.
This process involves an increasing number of Riemann solves as the dimensionality of the problem increases.
In the 1D case, CTU is identical to donor-cell upwind.


\section{Computing the Hybridization Coefficient}
We compute the hybridization coefficient $\eta$ using a modified multidimensional flux-corrected transport (FCT) algorithm.
Note that the time superscript notation ($n$) for fluxes is dropped from here out, but it is implied.
Our algorithm is based upon the method described first in \cite{Zalesak1979}.
Here is the generic FCT procedure:

\begin{enumerate}
\item Compute the high-order fluxes $\ang{\mb{F}^d_H}_{\pmhalf}$ over the cell volume $V_{\mb{i}}$.
\item Compute the low-order fluxes $\ang{\mb{F}^d_L}_{\pmhalf}$ and the corresponding low-order update $\ang{q}^{td}_{\mb{i}}$.
\begin{equation}
\ang{q}^{td}_{\mb{i}} = \ang{q}^{n}_{\mb{i}} - \frac{\Delta t}{h}\sum_d \left[ \ang{\mb{F}^d_L}_{\phalf} - \ang{\mb{F}^d_L}_{\mhalf} \right]
\end{equation}
\item Compute the antidiffusive fluxes $\ang{A^d}_{\pmhalf}$
\begin{equation}
\ang{A^d}_{\pmhalf} = \ang{\mb{F}^d_H}_{\pmhalf} - \ang{\mb{F}^d_L}_{\pmhalf}
\end{equation}
\item Limit the antidiffusive fluxes  
\begin{align} \label{lim_flux}
\ang{A^d_{\eta}}_{\pmhalf} &= \eta^d_{\pmhalf}\ang{A^d}_{\pmhalf} \\
&0\le \eta^d_{\pmhalf} \le 1 \nonumber
\end{align}
\item Update the solution with the limited antidiffusive fluxes
\begin{equation}
\ang{q}^{n+1}_{\mb{i}} = \ang{q}^{td}_{\mb{i}} - \frac{\Delta t}{h} \sum_d \left[ \ang{A^d_{\eta}}_{\phalf} - \ang{A^d_{\eta}}_{\mhalf} \right]
\end{equation}
\end{enumerate}

\subsection{Limiting the antidiffusive flux}
The primary challenge in the above formulation is comuting the hybridization coefficients ($\eta_{\pmhalf}$).
Following the procedure in \cite{Zalesak1979}, we compute the coefficients in the following manner: 

Pre-constrain the high-order fluxes $\ang{\mb{F}_H}_{\pmhalf}$.
This is a pre-limiting step that seeks to keep the high-order fluxes from generating new or steeper extrema.

Compute the sum ($P_{\mb{i}}^{\pm}$) of all the antidiffusive fluxes into and out of the cell and a measure of the diffusion ($Q_{\mb{i}}^{\pm}$).

\begin{align}
P_{\mb{i}}^{+} &= \sum_d^D \left[ \max \left( \ang{A^d}_{\mhalf},0 \right) - \min \left( \ang{A^d}_{\phalf},0 \right) \right] \\
Q_{\mb{i}}^{+} &= \left((q_{\max})_{\mb{i}} - \ang{q}_{\mb{i}}^{td}\right) \frac{h}{\Delta t} \\
P_{\mb{i}}^{-} &= \sum_d^D \left[ \max \left( \ang{A^d}_{\phalf},0 \right) - \min \left( \ang{A^d}_{\mhalf},0 \right) \right] \\
Q_{\mb{i}}^{-} &= \left(\ang{q}_{\mb{i}}^{td} - (q_{\min})_{\mb{i}}\right) \frac{h}{\Delta t}
\end{align}

Compute the least upper bounds $R_{\mb{i}}^{\pm}$.

\begin{align}
R_{\mb{i}}^{+} &=
  \begin{cases}
	\min \left( 1.0 , Q_{\mb{i}}^{+} / P_{\mb{i}}^{+} \right) & \text{if } P_{\mb{i}}^{+} > 0.0 \\
	0.0 & \text{otherwise}
  \end{cases} \\
R_{\mb{i}}^{-} &=
  \begin{cases}
	\min \left( 1.0 , Q_{\mb{i}}^{-} / P_{\mb{i}}^{-} \right) & \text{if } P_{\mb{i}}^{-} > 0.0 \\
	0.0 & \text{otherwise}
  \end{cases} 
\end{align}

Select the hybridization coefficient with the most restrictive upper bound.

\begin{align}
\eta_{\phalf} &=
  \begin{cases}
   \min \left( R_{\splus}^{+} , R_{\mb{i}}^{-} \right)  & \text{if }  \ang{A^d}_{\phalf} > 0.0 \\
   \min \left( R_{\mb{i}}^{+} , R_{\splus}^{-} \right)  & \text{if }  \ang{A^d}_{\phalf} \le 0.0
  \end{cases}
\end{align}

In the above description the user is provided with two design choices: pre-constraint for the high-order flux and method of computing the solution bounds 
$(q_{\max})_{\mb{i}}$ and $(q_{\min})_{\mb{i}}$.

\subsection{Computing the Solution Bounds}
Compute initial estimates of the solution bounds, $(q_{\max})_{\mb{i}}$ and $(q_{\min})_{\mb{i}}$.
First, compute the bounded solutions in a rectangular stencil ($B_{\mb{i}}$) that is $[2s_{\mb{i}}+1]^D$ cells in size, where $s_{\mb{i}}$ is the stencil size. In this study, $s_{\mb{i}}$ was allowed to vary depending on the local velocity field in the following manner:

\begin{align}
s_{\mb{i}} &=
  \begin{cases}
   2  &\text{if } \ \ \sigma \max_{d} \limits |(\mb{u}^d)_{\mb{i}}| \ \ge \ 0.5 \\
   1  &\text{otherwise }  
  \end{cases}
\end{align}

where $\sigma = \max (|\vec{\mb{u}}|) \Delta t/h$ is the global CFL number.

The stencil was allowed to vary to ensure good performance over the spectrum of possible CFL numbers.
It was found that at low CFL numbers, $s_{\mb{i}} = 2$ lead to excess diffusion near discontinuities.
Alternatively at high CFL numbers, $s_{\mb{i}} = 1$ lead to excess diffusion.
Allowing the stencil size to vary mitigated these issues.

After the stencil is determined, four bounds are computed: (1) max based on $\ang{q}^n$, (2) min based on $\ang{q}^n$, (3) max based on $\ang{q}^{td}$, and (4) min based on $\ang{q}^{td}$.

\begin{align}
 &(q_{\max})_{\mb{i}}^{n} =  \max \left( B_{\mb{i}} (\ang{q}^n) \right) \\
 &(q_{\min})_{\mb{i}}^{n} =  \min  \left( B_{\mb{i}} (\ang{q}^n)  \right) \\
 &(q_{\max})_{\mb{i}}^{td} =  \max \left( B_{\mb{i}} (\ang{q}^{td}) \right) \\
 &(q_{\min})_{\mb{i}}^{td} =  \min \left(  B_{\mb{i}} (\ang{q}^{td}) \right)
 \end{align}
 
 Then select the upper and lower bound of the two estimates.
 
 \begin{align}
 &(q_{\max})_{\mb{i}} =  \max( (q_{\max})_{\mb{i}}^{n},  (q_{\max})_{\mb{i}}^{td} ) \\
 &(q_{\min})_{\mb{i}} =  \min( (q_{\min})_{\mb{i}}^{n} ,  (q_{\min})_{\mb{i}}^{td} )
\end{align}

\subsubsection{Accurate Solution Bounds at Smooth Extremum}
For the vast majority of cells within the domain, the previous bound computation is sufficiently accurate.
However, computing bounds at extremum is more complicated.
Ideally the bounds need to keep the solution monotonic and positive, but the bounds should also not ``clip'' the solution.
There are few different methods for avoiding clipping, and we use a geometric construction that is only applied at smoothly varying extrema.
It is based on the ideas in \cite{Colella2008}.

The smooth extremum criteria in 1D is:

\begin{align}
(ext^d)_{\mb{i}} \ = \ &\min [ \ (dq)_{\mb{i}}\cdot (dq)_{\mb{i}+\mb{e}^d} \ , \ (dq)_{\mb{i}-\mb{e}^d}\cdot (dq)_{\mb{i}+2\mb{e}^d} \ ] \le 0.0  \ \ \&\& \\ \nonumber
& 1.25\cdot (dqtot)_{\mb{i}} < (tv)_{\mb{i}}
\end{align}

where the following definitions are used:

\begin{align}
(dq)_{\mb{i}} &= \ang{q}_{\mb{i}}^{td} - \ang{q}_{\mb{i}-\mb{e}^d}^{td} \\
(dqtot)_{\mb{i}} &= |\ang{q}_{\mb{i}+2\mb{e}^d}^{td} - \ang{q}_{\mb{i}-2\mb{e}^d}^{td} | \\
(tv)_{\mb{i}} &= |(dq)_{\mb{i}+2\mb{e}^d}| + |(dq)_{\mb{i}+\mb{e}^d}| + |(dq)_{\mb{i}}| + |(dq)_{\mb{i}-\mb{e}^d}| 
\end{align}

This criteria has two parts.
First, check for a sign change in the first derivative.
The sign change will indicate either an extremum or a discontinuity in the solution.
Second, ensure that the solution locally is not a perterbation of a discontinuity.

For a smooth multidimensional extremum either $(ext^d)_{\mb{i}}$ must be true in all dimensions or it must be true for some $d$ and the solution must remain constant along the dimensions in which $(ext^d)_{\mb{i}}$ is not true.
We use the following criterion to determine if the solution is constant:

\begin{equation}
\max(|(q^d_{\max})_{\mb{i}}^{td} - \ang{q}_{\mb{i}}^{td}|,|(q^d_{\min})_{\mb{i}}^{td} - \ang{q}_{\mb{i}}^{td}|) \le 10^{-14} 
\end{equation}

where the following definitions hold

\begin{align}
(q^d_{\max})_{\mb{i}}^{td} &= \max \left ( \ang{q}^{td}_{\sminus} , \ang{q}^{td}_{\mb{i}}, \ang{q}^{td}_{\splus} \right )  \\
(q^d_{\min})_{\mb{i}}^{td} &= \min \left (   \ang{q}^{td}_{\sminus} , \ang{q}^{td}_{\mb{i}}, \ang{q}^{td}_{\splus}  \right )
\end{align}

Once we have determined that the solution at $V_{\mb{i}}$ is at a smooth extremum, we compute new values of $(q_{\max})_{\mb{i}}$ and $(q_{\min})_{\mb{i}}$.
The first step is to construct a parabolic function from the local values of $\ang{q}$ centered at $\mb{x}_{\mb{i}}$:

\begin{equation} \label{quad}
q^d(x) =  \left( \frac{(d2q)_{\mb{i}}^n}{2} \right) x^2 + \left( \frac{( \ang{q}_{\splus} - \ang{q}_{\sminus} )}{2} \right) x + \ang{q}_{\mb{i}}
\end{equation}

where

\begin{equation}
(d2q)_{\mb{i}}^n = \ang{q}_{\splus}^n + \ang{q}_{\sminus}^n - 2\ang{q}_{\mb{i}}^n
\end{equation}

The location of the vertex ($x_c$) is given by the ratio $-b/2a$, where a and b are the quadratic and linear coefficients (\ref{quad}):

\begin{equation}
x_c = -\frac{\ang{q}_{\splus} - \ang{q}_{\sminus} }{2(d2q)_{\mb{i}}^n}
\end{equation}

and $-0.5 \le x_c \le 0.5$.
Then, evaluate the quadratic at the vertex to find the extremum value as well as deconvolve to get an estimate of the point value:

\begin{equation}
(q^d_{ext})_{\mb{i}} = \left( \frac{(d2q)_{\mb{i}}^n}{2} \right) x_c^2 + \left( \frac{( \ang{q}_{\splus} - \ang{q}_{\sminus} )}{2} \right) x_c + \ang{q}_{\mb{i}} - \frac{(d2q)_{\mb{i}}^n}{24}
\end{equation}

Select the largest $(q^d_{ext})_{\mb{i}}$ or smallest $(q^d_{ext})_{\mb{i}}$ depending on the sign of the second derivative:

\begin{align}
(q_{ext})_{\mb{i}} &=
  \begin{cases}
   \max_d \limits ((q^d_{ext})_{\mb{i}} , (q_{\max})_{\mb{i}})  &\text{if } \ \ \sgn((d2q)_{\mb{i}}^n) \le 0.0 \\
   \min_d \limits  ((q^d_{ext})_{\mb{i}} , (q_{\min})_{\mb{i}}) &\text{otherwise }  
  \end{cases}
\end{align}

Finally, compute the appropriate extremum bound by augmenting the solution value at the previous time by a scaled difference between the extremum value and the solution value:

\begin{align}
(q_{\max})_{\mb{i}} &=
  \begin{cases}
    q_{\mb{i}}^n + \max(0.0 , 2.0 \left | (q_{ext})_{\mb{i}} - q_{\mb{i}}^n \right | )  &\text{if } \ \ \sgn((d2q)_{\mb{i}}^n) \le 0.0 \\
    (q_{\max})_{\mb{i}} &\text{otherwise }  
  \end{cases} \\
(q_{\min})_{\mb{i}} &=
  \begin{cases}
   q_{\mb{i}}^n + \min(0.0 , 2.0 \left | (q_{ext})_{\mb{i}} - q_{\mb{i}}^n \right | )  &\text{if } \ \ \sgn((d2q)_{\mb{i}}^n) > 0.0 \\
   (q_{\min})_{\mb{i}} &\text{otherwise }  
  \end{cases}
\end{align}

\subsubsection{Updating $R_{\mb{i}}^{\pm}$ at Extrema}
We flag the extrema at which the Laplacian is changing sign. We compute the $D$ dimensional approximation to the Laplacian over a 3-point stencil:

\begin{equation}
\Delta q_{\mb{i}} = \sum^D_{d = 1} {\frac{\partial^2 q_{\mb{i}}}{\partial x_d^2}} \approx \sum^D_{d = 1} \frac{(d2q)_{\mb{i}}^n}{h^2}
\end{equation}

If $\Delta q$ changes sign anywhere in the 3-point vicinity of $\mb{i}$, then we flag that cell $\mb{i}$.
We then update the least upper bound multiplier at the flagged cells:

\begin{equation}
R_{\mb{i}}^{\pm} = 0 \ \ \text{if } \mb{i} \text{ flagged} 
\end{equation}

\subsection{Pre-constraining the High-Order Flux}

We pre-constrain the high-order flux by modifying the high-order fluxes where they would otherwise accentuate or produce a new extremum.
In practice, the value of the antidiffusive flux is edited instead of the high-order flux directly.
Following \cite{Zalesak1979} we set the antidiffusive flux to zero in these regions.

The baseline condition for applying the pre-constraint:

\begin{equation}
\ang{A^d}_{\phalf}(\ang{q}^{td}_{\splus} - \ang{q}^{td}_{\mb{i}}) \le 0.0
\end{equation}

However, this condition alone is not sufficient.
This condition will occasionally be satisfied at smooth areas in the solution.
We add the following conditions to make sure we only apply this condition away from smooth areas:

\begin{equation}
\min \left[ (d2q)_{\splus}^n\cdot (d2q)_{\mb{i}}^n \ , \ (d2q)_{\mb{i}}^n\cdot (d2q)_{\sminus}^n \ , \ (d2q)_{\splus}^n\cdot (d2q)_{\mb{i}+2\mb{e}^d}^n \right] < 0.0
\end{equation}
\begin{equation} \label{smooth_preconstraint}
|\ang{A^d}_{\phalf}| \le \frac{|(\mb{u}^d)_{\phalf}| h}{2} \ \left( 1-\sigma_{\phalf} \right) \ \frac{\ |(d2q)_{\mb{i}} + (d2q)_{\splus}|}{2} 
\end{equation}

where $\sigma_{\phalf} = |(\mb{u}^d)_{\phalf}| \Delta t / h$.

The first constraint above attempts to detect a discontinuity in the solution.
Discontinuities are natural places for high-order fluxes and consequently antidiffusive fluxes to generate a new extremum.
However, there are smooth multidimensional solutions in which the second derivative naturally changes sign.
The second constraint seeks to preclude this case.
The term on the right hand side of the inequality (\ref{smooth_preconstraint}) is the d-directional dissipation term, scaled by the cell size, in the modified equation analysis (\ref{mod_analysis}) of CTU applied to the advection equation 

\begin{equation} \label{mod_analysis}
\pd{q}{t} + \sum_d^D \left( \mb{u}^d \pd{q}{x_d} \right) = \sum_d^D \left( \frac{\mb{u}^d h}{2} (1-\sigma_d) \frac{\partial^2 q}{\partial x_d^2} \right) + O(h^2)
\end{equation}

We are interested in comparing the magnitude of this dissipative term to the antidiffusive flux.
While it is true that the second derivative takes on a zero value right at a discontinuity, the second derivative takes on large values immeadiately surrounding this point.
This means that the dissipative term is large in the neighborhood of discontinuities.
If the magnitude of the disspative term is large relative to the the antidiffusive term, then we assume we are near a discontinuity and allow the pre-constraint.

\section{Results}

Results in one and two dimensions are presented.
A total of four initial conditions were investigated.
Of the four, one initial condition was smooth and the others contained a discontinuity.
For the two dimensional tests, we used two different velocity fields (\ref{v_field_1}-\ref{v_field_2}): constant diagonal and solid body rotation.

\begin{align} \label{v_field_1}
\mb{u} &= [1,1] \\ \label{v_field_2}
\mb{u} &= 2\pi [y - 0.5, 0.5 - x]
\end{align}

The center for the constant velocity initial condition was in the middle of the domain, whereas it was offset by 0.25 of the grid height for the solid body rotation examples

\begin{align*}
\mb{x}^{\text{const}}_{c} &= \left( 0.5 , 0.5 \right) \\
\mb{x}^{\text{solid}}_{c} &= \left( 0.5 , 0.75 \right)
\end{align*}

\subsection{Initial Conditions}
The smooth initial condition was constructed as a power of cosines.

\begin{equation}
  q_{\mb{i}}(t_0) =
  \begin{cases}
   \cos^8 \left( \frac{\pi}{2} \cdot \frac{R}{R_0} \right) & \text{if }  R \le R_0 \\
   0 & \text{otherwise}
  \end{cases}
\end{equation}
\begin{align*}
\mb{x}_{\mb{i}} &\in [0,1] \\
R &= \sqrt{(\mb{x}_{\mb{i}} - \mb{x}_{c})^2}
\end{align*}

For this initial condition, $R_0 = 15$ was chosen.

Three different discontinuous initial conditions were investigated.
The first was a square and is described as:

\begin{align}
  \label{square_ic}
  q_{\mb{i}}(t_0) &=
  \begin{cases}
   1 & \text{if  for each $D$: }  |x^D_{\mb{i}} - x^D_{c}| \le 0.15 \\
   0 & \text{otherwise}
  \end{cases}
\end{align}

The next is a semi-ellipse:

\begin{align}
  \label{semi_ic}
  q_{\mb{i}}(t_0) &=
  \begin{cases}
   \sqrt{1.0 - \left(\frac{R}{R_0}\right)^2} & \text{if  }  R \le R_0 \\
   0 & \text{otherwise}
  \end{cases}
\end{align}

In this case $R_0 = 0.25$.

The last test case is the classic slotted cylinder in two dimensions \cite{Zalesak1979}.
For this paper we used $N = 256$ to represent the cylinder.
This is roughly twice as resolved as the original.

\subsection{One-Dimensional Tests}
\begin{figure}[tb]
\centering
	\begin{subfigure}[b]{\sfwidth}
		\resizebox{\textwidth}{!}{\input{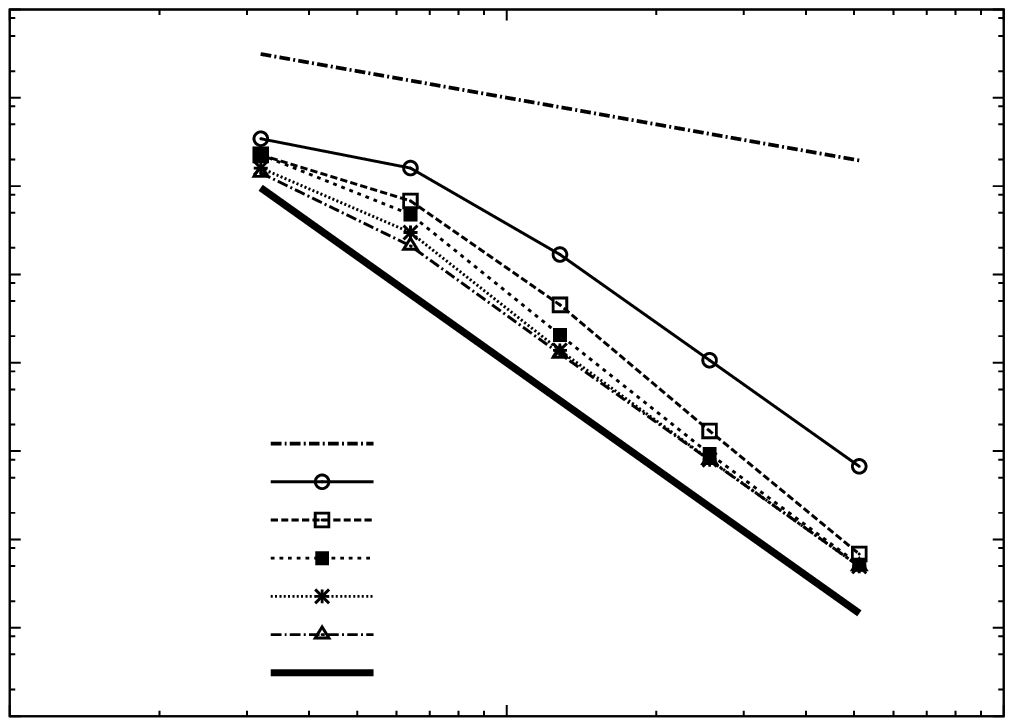}}
		\caption{Solution Error without Limiter}
		\label{fig:convergence}
	\end{subfigure}
	\begin{subfigure}[b]{\sfwidth}
		\resizebox{\textwidth}{!}{\input{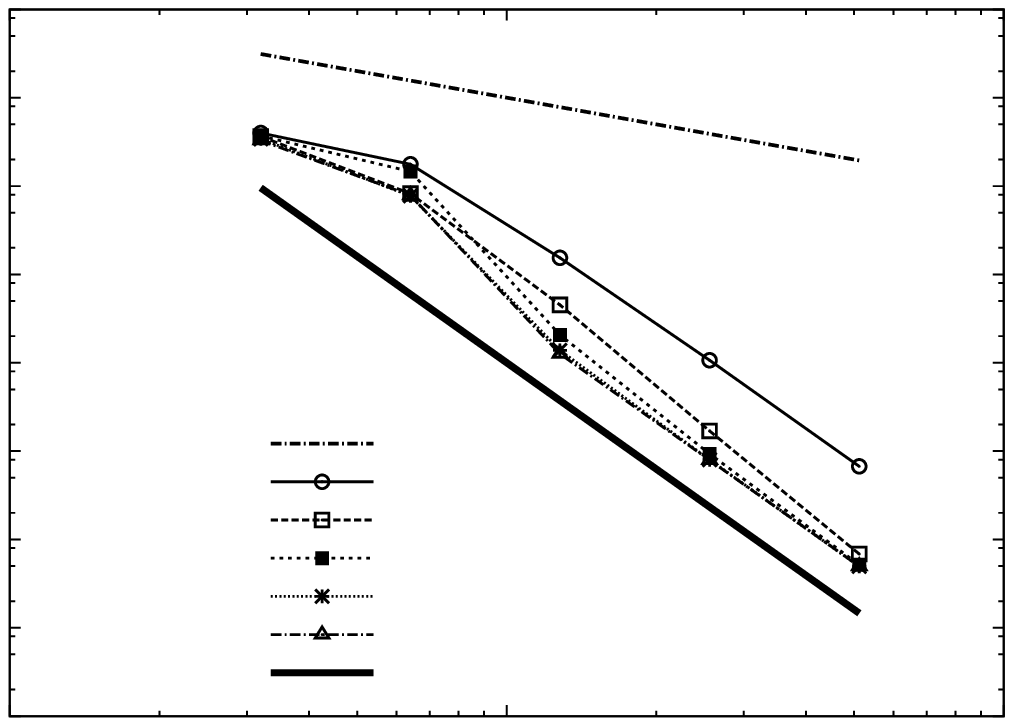}}
		\caption{Solution Error with Limiter}
		\label{fig:solution}
	\end{subfigure}
\captionsetup{width=\textwidth}
\caption[convergence_1d]{Max Norm Errors in 1D}
\label{convergence_1d}
\end{figure}

The first requirement for the limiter method is that it reduces to the high-order scheme away from discontinuities.
All of the high-order schemes achieved similar errors for smooth solutions (\ref{convergence_1d}).
The rate of convergence for each method was $4.0$, which is the expected rate due to the usage of RK4 for the time integration scheme.
This implies that the limiter is not being activated in smooth regions, and the first requirement is met.

\begin{figure}[tb]
\centering
	\begin{subfigure}[b]{\sfwidth}
		\resizebox{\textwidth}{!}{\input{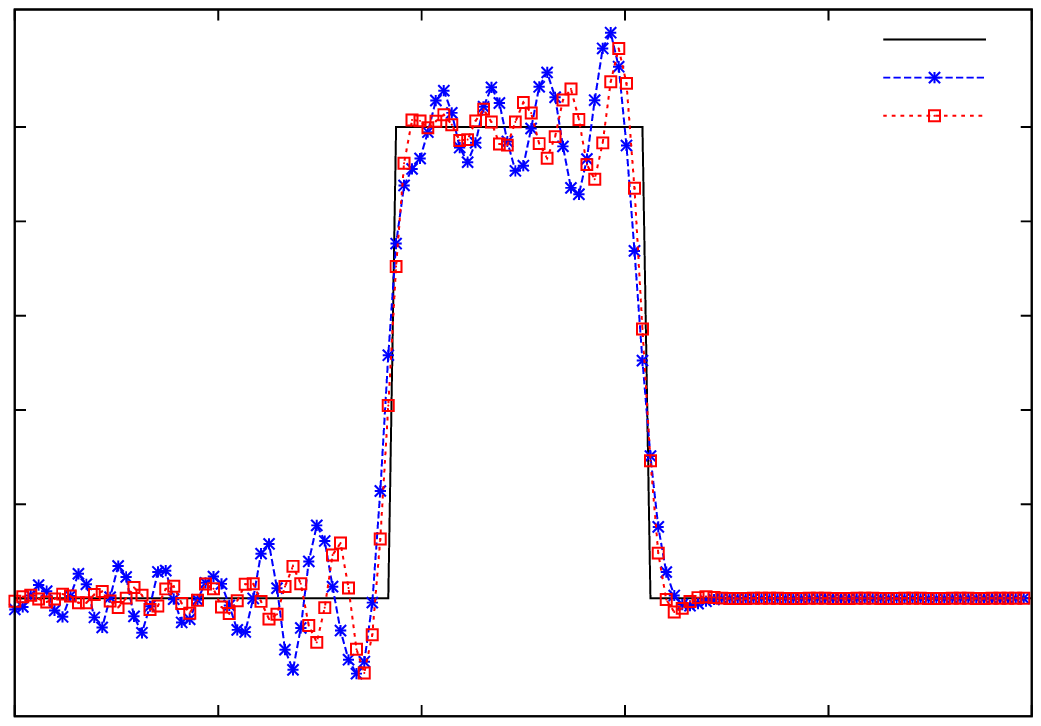}}
		\caption{Centered Differences}
		\label{fig:square.11}
	\end{subfigure}
	\begin{subfigure}[b]{\sfwidth}
		\resizebox{\textwidth}{!}{\input{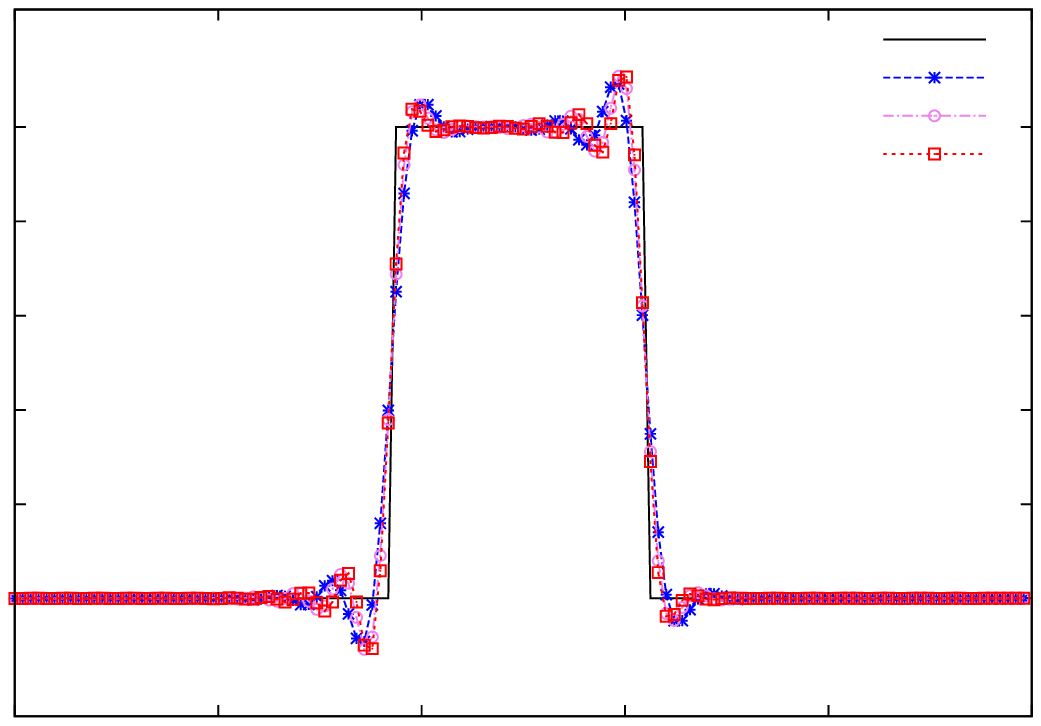}}
		\caption{Upwind Differences}
		\label{fig:square.12}
	\end{subfigure}
\captionsetup{width=\textwidth}
\caption[square_1nolim]{Comparison of Square Solutions with Different High-Order Spatial Differencing Schemes ($\sigma=0.8$, $t=1.0$, $N=128$)}
\label{square_1nolim}
\end{figure}

The second requirement is that the limiter method accurately represent discontinuities.
To determine if this requirement is met, we analyzed the performance of the high-order schemes with and without the limiter.
The first test was a comparison between the high-order schemes without the limiter (\ref{square_1nolim}).
In the presence of a discontinuity, the centered difference solutions produced more pronounced oscillations than the upwind solutions.
This result is consistent with the amplitude and phase error analysis presented earlier.

\begin{figure}[tb]
\centering
	\begin{subfigure}[b]{\sfwidth}
		\resizebox{\textwidth}{!}{\input{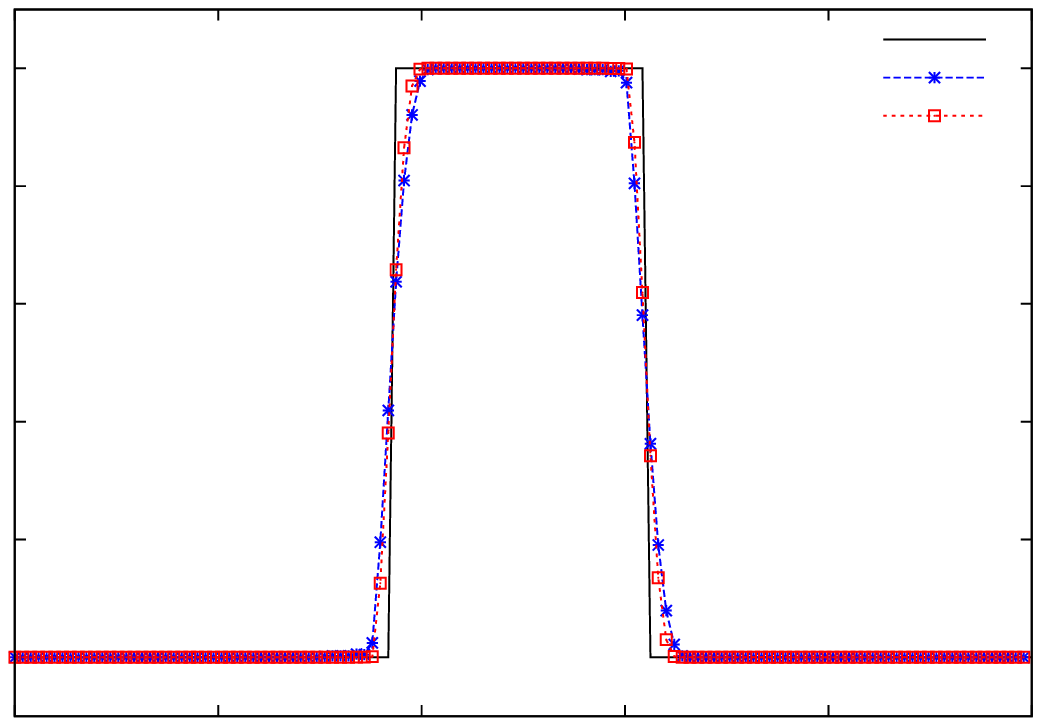}}
		\caption{Centered Differences}
		\label{fig:square.13}
	\end{subfigure}
	\begin{subfigure}[b]{\sfwidth} 
		\resizebox{\textwidth}{!}{\input{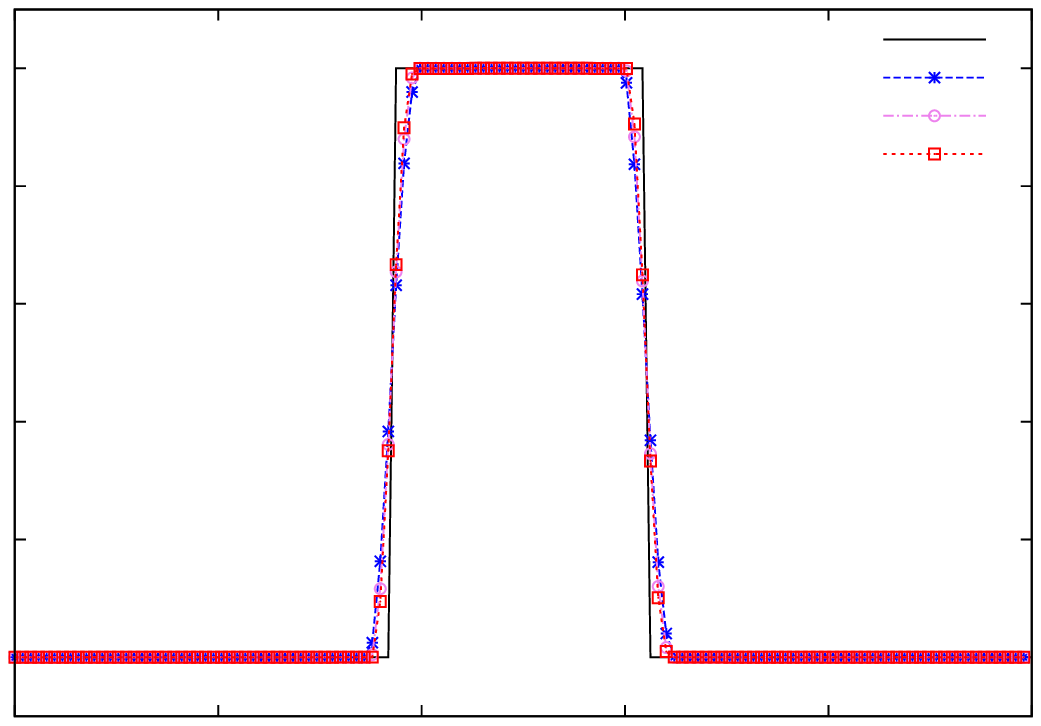}}
		\caption{Upwind Differences}
		\label{fig:square.14}
	\end{subfigure}
\captionsetup{width=\textwidth}
\caption[square_1nolim]{Comparison of Limited Square Solutions with Different High-Order Spatial Differencing Schemes ($\sigma=0.8$, $t=1.0$, $N=128$)}
\label{square_1lim}
\end{figure}

\begin{figure}[tb]
\centering
	\begin{subfigure}[b]{\sfwidth}
		\resizebox{\textwidth}{!}{\input{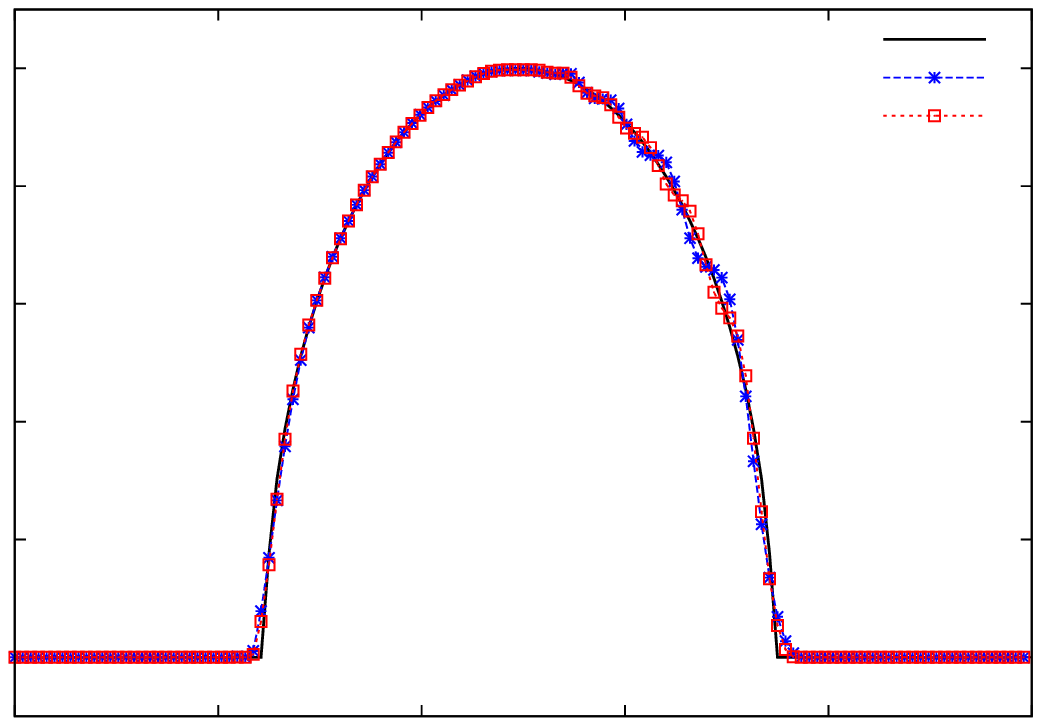}}
		\caption{Centered Differences}
		\label{semi.11}
	\end{subfigure}
	\begin{subfigure}[b]{\sfwidth}
		\resizebox{\textwidth}{!}{\input{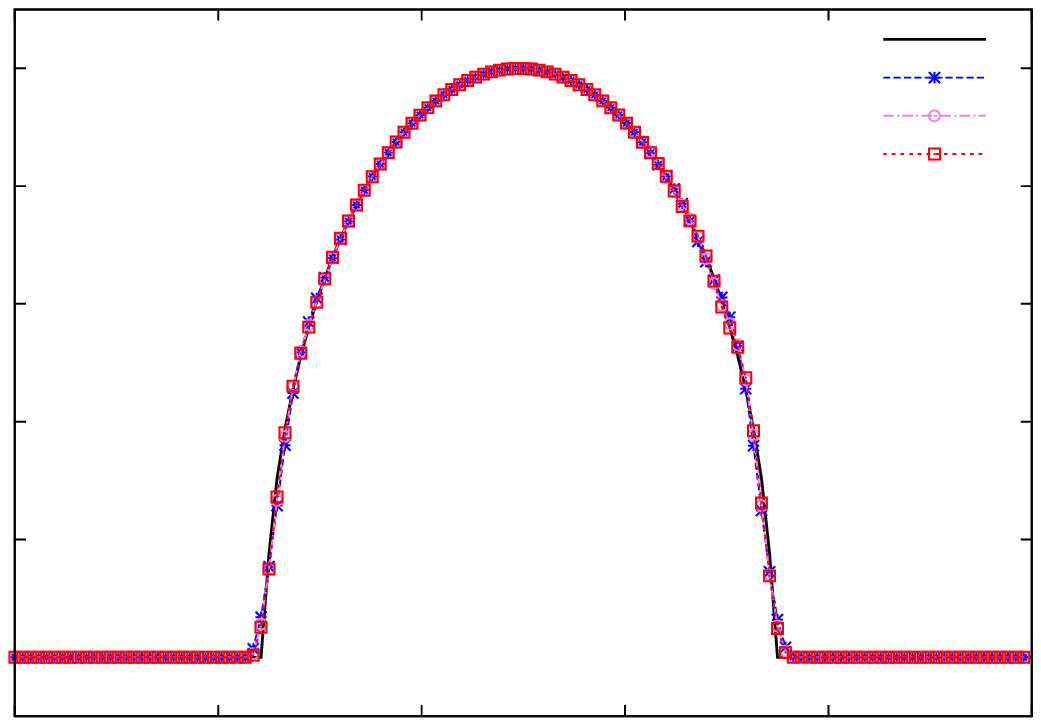}}
		\caption{Upwind Differences}
		\label{semi.12} 
	\end{subfigure}
\captionsetup{width=\textwidth}
\caption[semi_1lim]{Comparison of Limited Semi-Ellipse Solutions with Different High-Order Spatial Differencing Schemes ($\sigma=0.8$, $t=1.0$, $N=128$)}
\label{semi_1lim}
\end{figure}

The more important comparisons are how these high-order schemes perform with the limiter on.
For the square initial condition, there is very little difference in the solution (\ref{square_1lim}).
The only noticeable difference is some additional diffusion on the tailing edge of the centered difference solution.
However, for the semi-ellipse, the centered difference schemes produce oscillations even in the limited solution near the discontinuity (\ref{semi_1lim}).
The oscillations are bounded by FCT; however, they are clearly undesirable.
The upwind schemes also produce some oscillations, but again the magnitude of the oscillations is smaller.
This is motivation for selecting the upwind methods moving forward.

\subsection{Two Dimensional Tests}
\begin{figure} [tb]
\centering
	\begin{subfigure}[b]{\sfwidth}
		\resizebox{\textwidth}{!}{\input{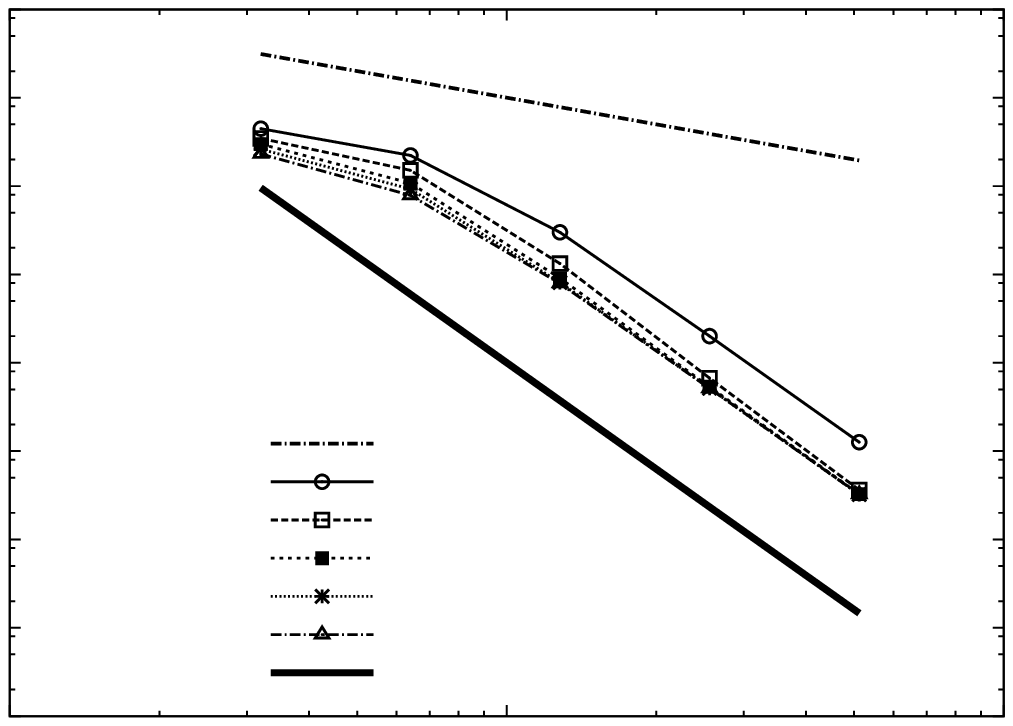}}
		\caption{Constant Velocity without Limiter}
		\label{error_const_2d}
	\end{subfigure}
	\begin{subfigure}[b]{\sfwidth}
		\resizebox{\textwidth}{!}{\input{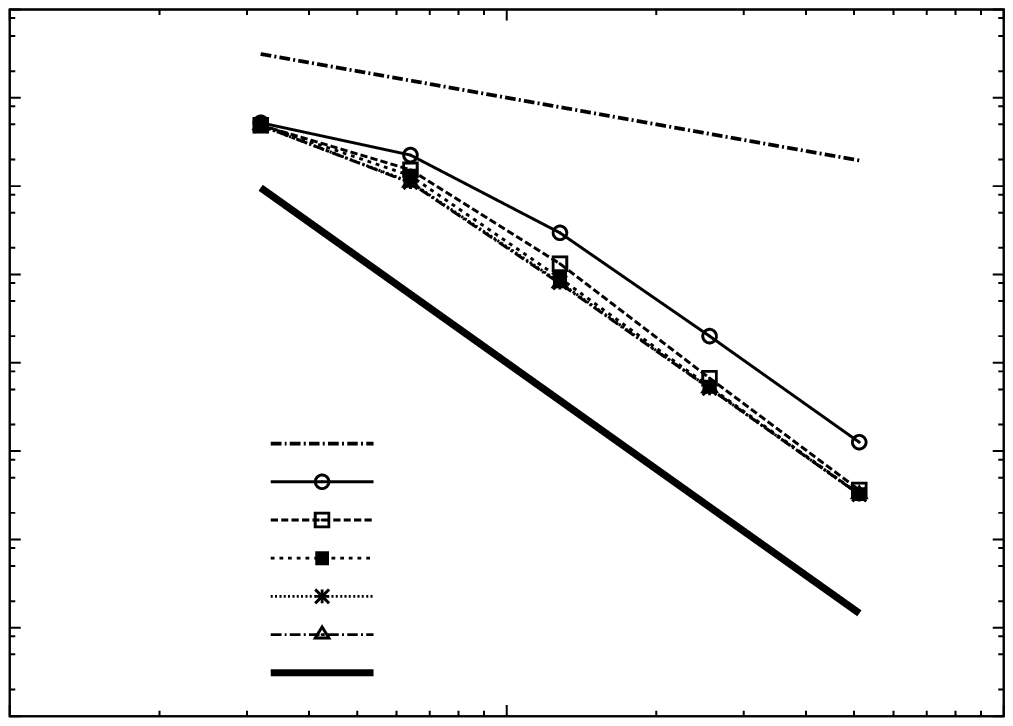}}
		\caption{Constant Velocity with Limiter}
		\label{error_const_2d_lim}
	\end{subfigure} \\
	\begin{subfigure}[b]{\sfwidth}
		\resizebox{\textwidth}{!}{\input{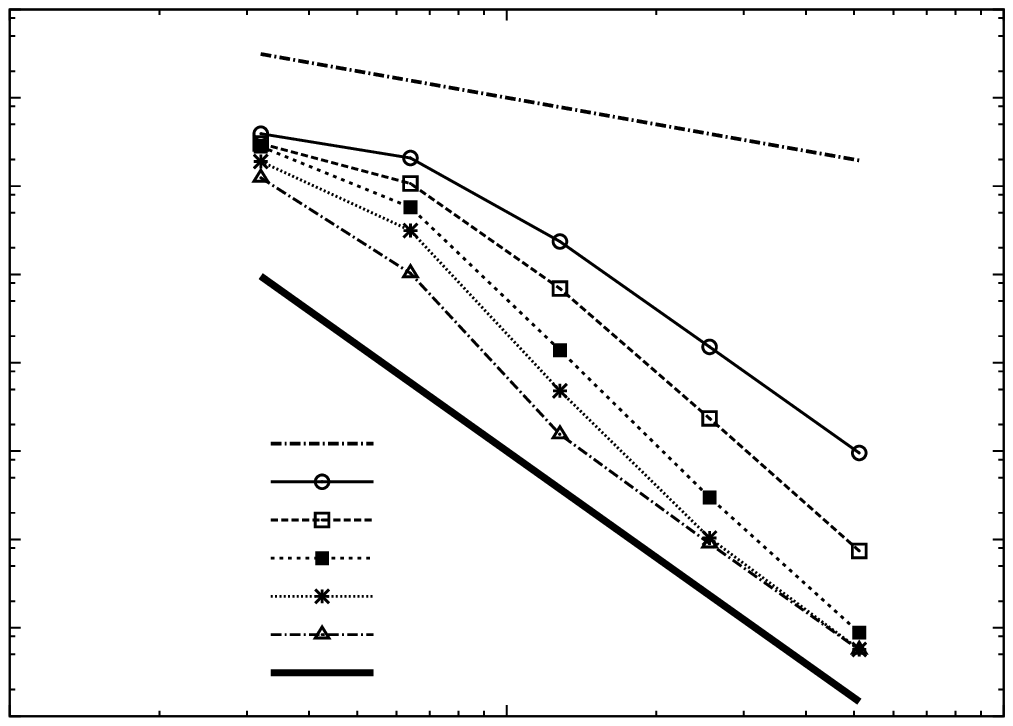}}
		\caption{Rotation Velocity without Limiter}
		\label{error_rot_2d}
	\end{subfigure}
	\begin{subfigure}[b]{\sfwidth}
		\resizebox{\textwidth}{!}{\input{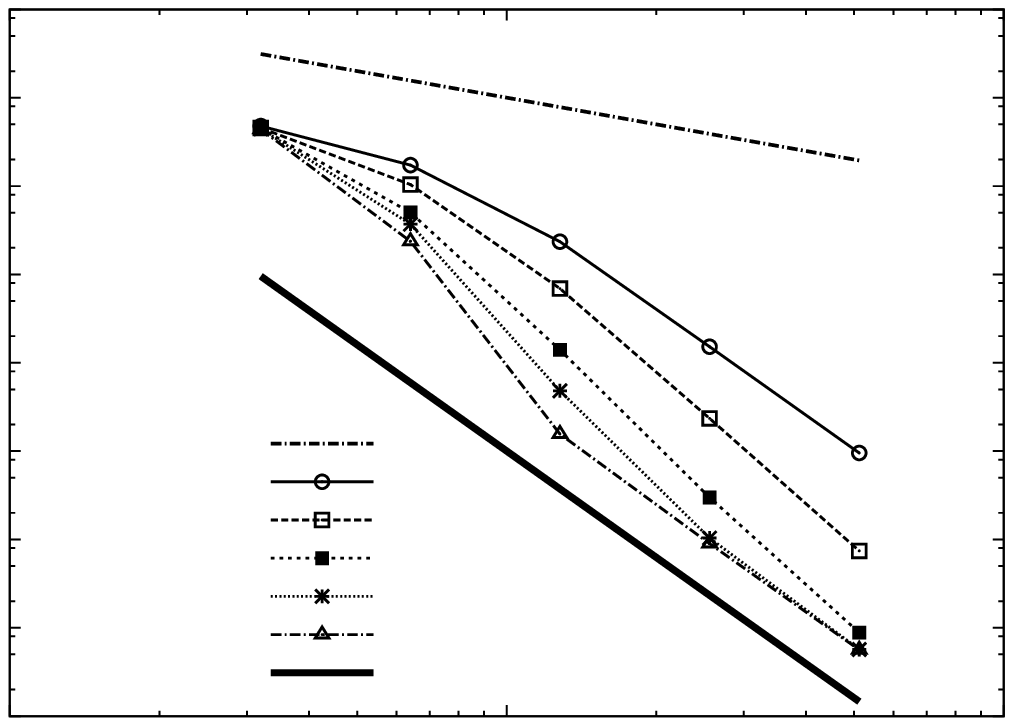}}
		\caption{Rotation Velocity with Limiter}
		\label{error_rot_2d_lim}
	\end{subfigure} 
\caption[error_2d]{Max Norm Errors in 2D}
\label{error_2d}
\end{figure}

The smooth test errors in two dimensions are reported with the limiter for both velocity fields (\ref{error_2d}).
As in the one dimensional case, the high-order solution accuracy requirement is met.
Each method achieves fourth order accuracy, at a minimum.
Interestingly, for the solid body rotation solution, the error reduction is greater than fourth order.
Since the spatial differencing error is sixth order accurate, this suggests that the time discretization error is smaller than expected even for large values of $\sigma$.
Also the ninth order scheme is still convergent right at its theoretical stability limit ($\sigma \sim 0.8$).

\begin{figure} [tb]
\centering
	\begin{subfigure}[b]{\sfwidth}
		\inputImgEps{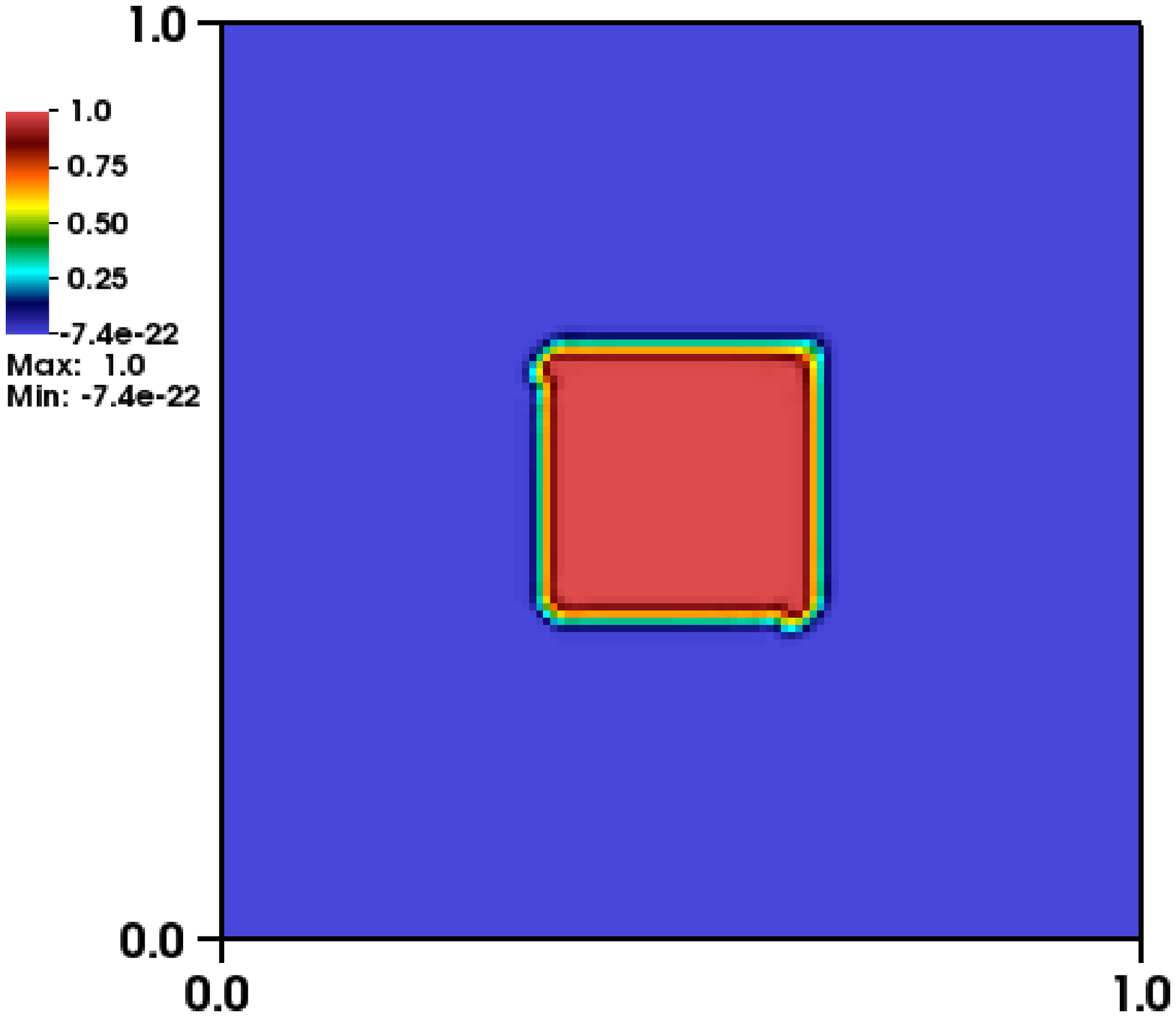}
		\caption{Constant Velocity Solution}
		\label{square_const_22}
	\end{subfigure}
	\begin{subfigure}[b]{\sfwidth}
		\inputImgEps{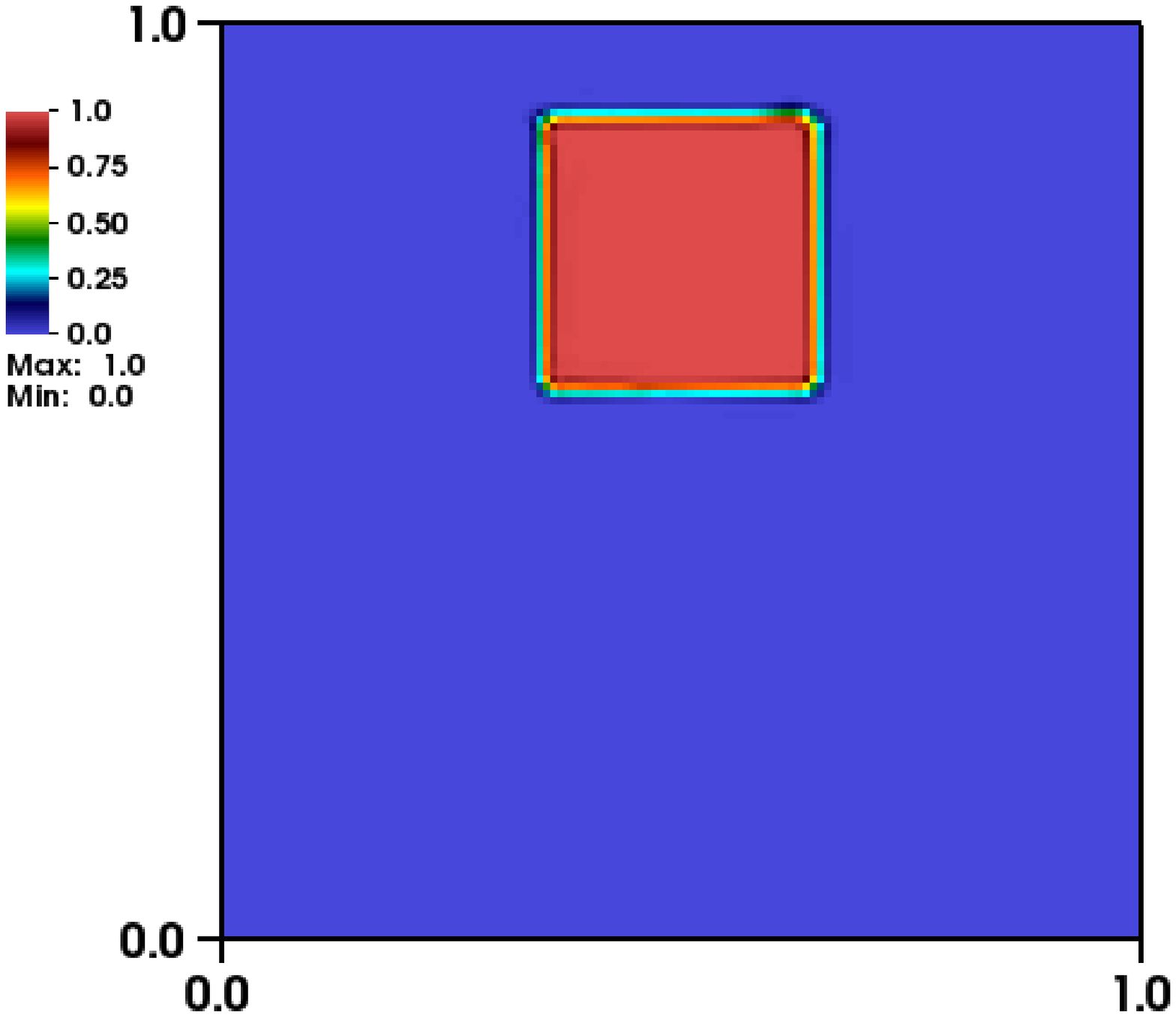}
		\caption{Solid Body Rotation Solution}
		\label{square_rot_22}
	\end{subfigure} \\
	\begin{subfigure}[b]{\sfwidth}
		\resizebox{\textwidth}{!}{\input{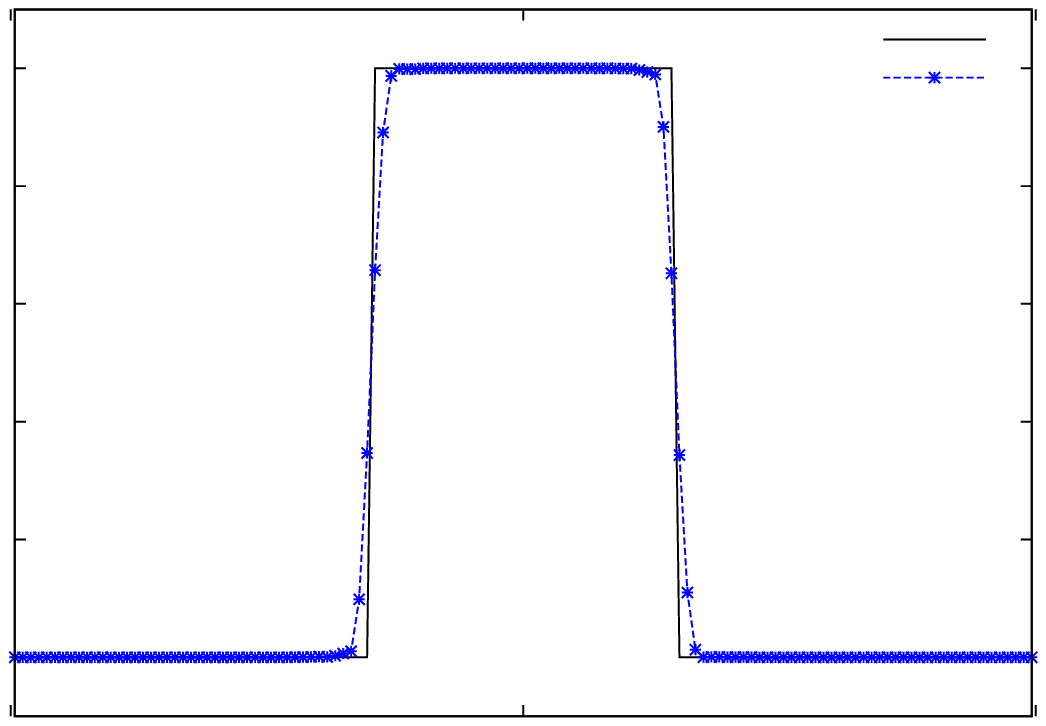}}
		\caption{Centerline Comparison}
		\label{square_const_23}
	\end{subfigure}
	\begin{subfigure}[b]{\sfwidth}
		\resizebox{\textwidth}{!}{\input{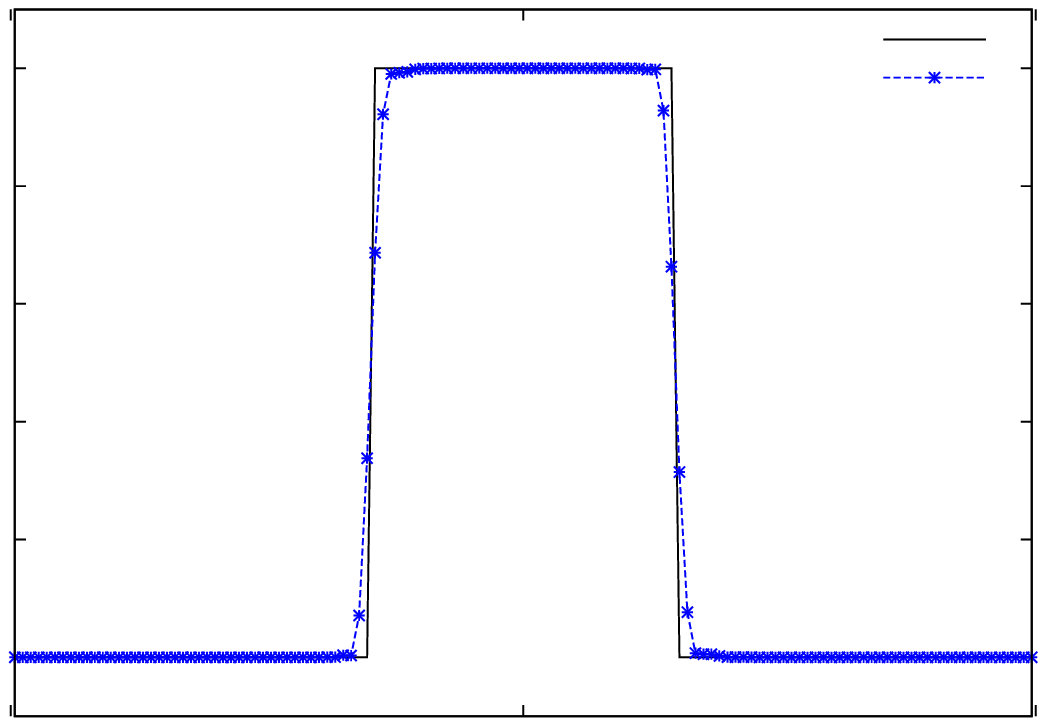}}
		\caption{Centerline Comparison}
		\label{square_rot_23}
	\end{subfigure}
\caption[square_2d]{Square Solutions ($t=1.0,\ N=128$, $\sigma = 0.8$)}
\label{square_2d}
\end{figure}

\begin{figure} [tb]
\centering
	\begin{subfigure}[b]{\sfwidth}
		\inputImgEps{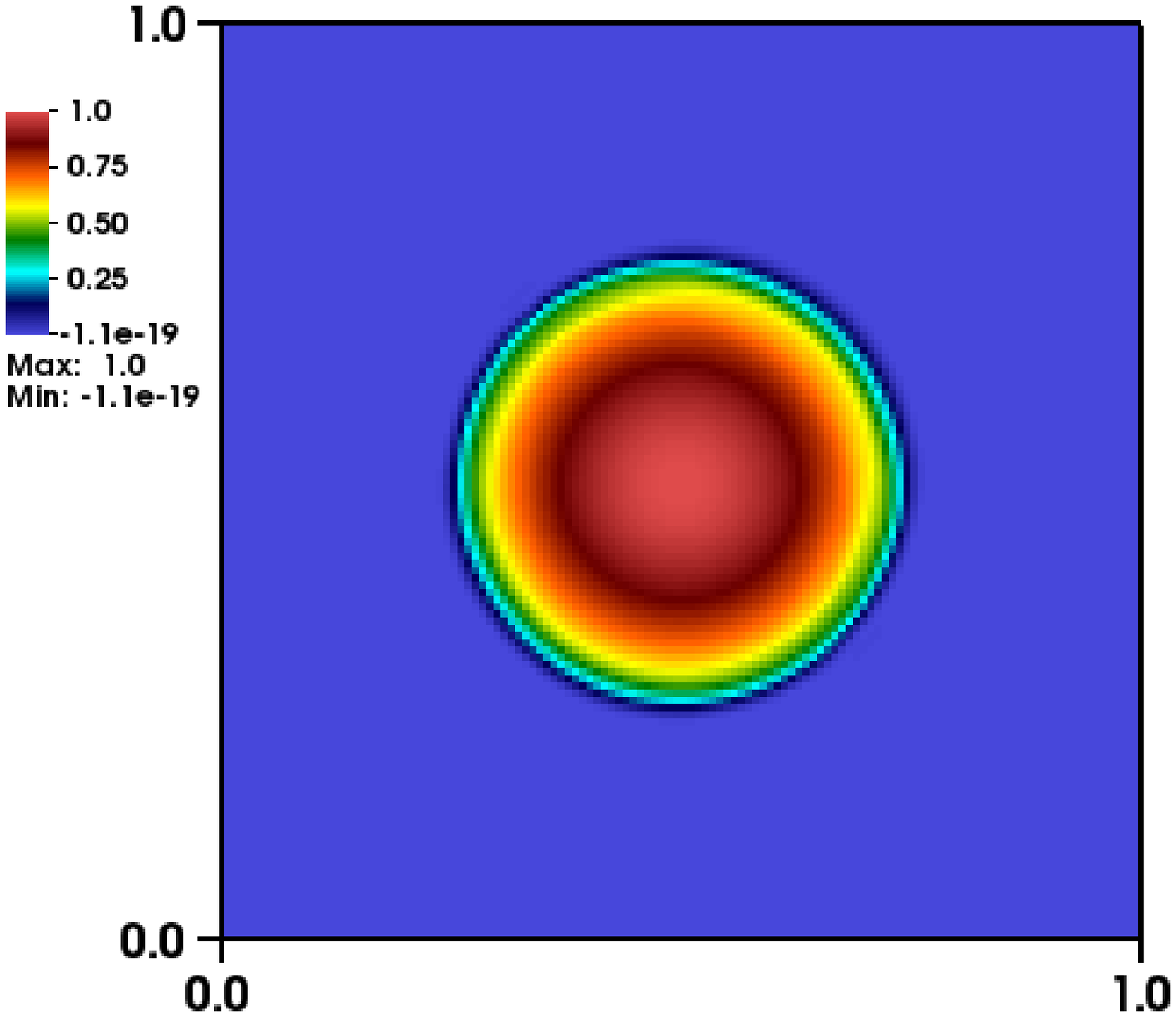}
		\caption{Constant Velocity Solution}
		\label{semi_const_22}
	\end{subfigure}
	\begin{subfigure}[b]{\sfwidth}
		\inputImgEps{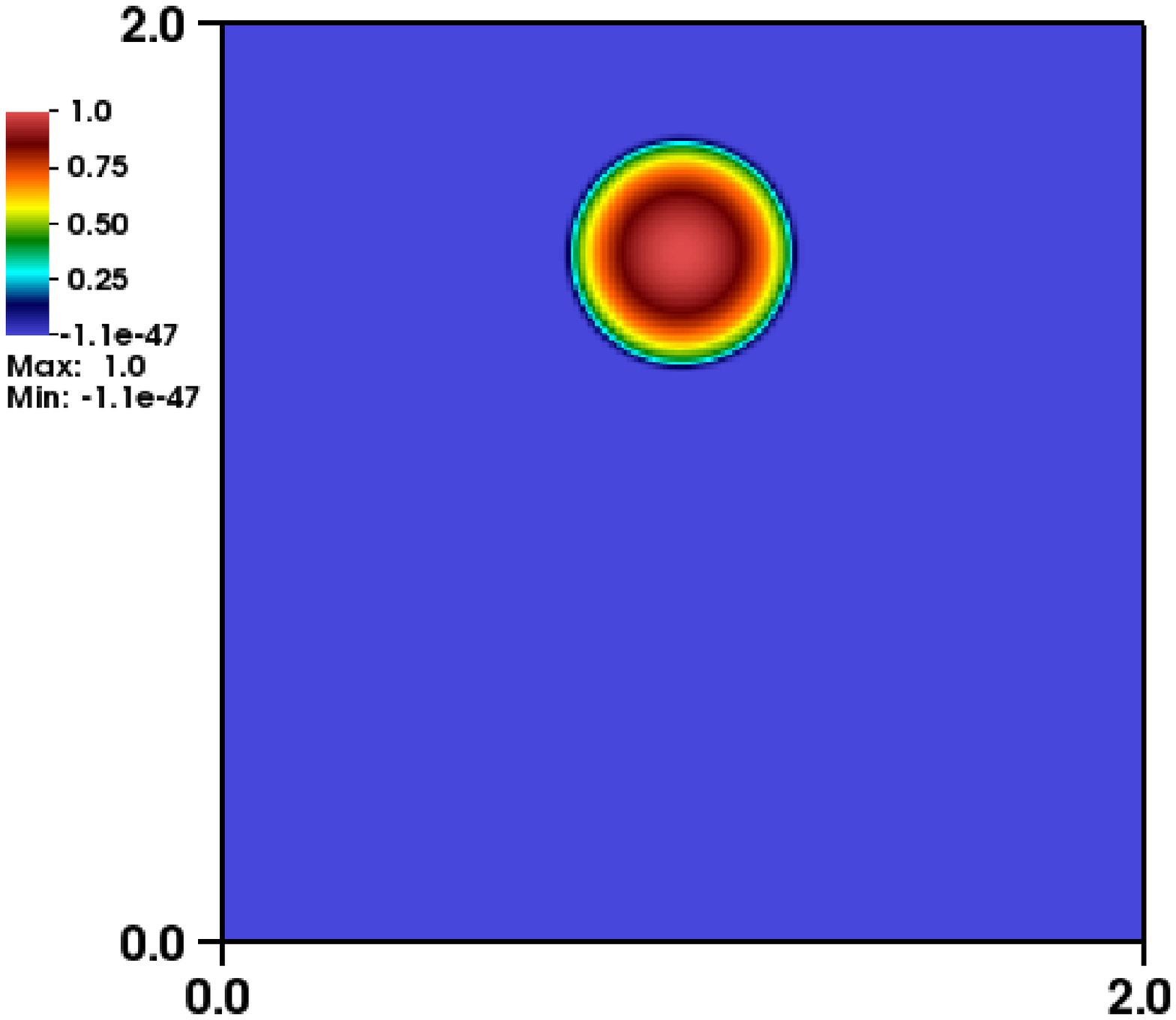}
		\caption{Solid Body Rotation Solution}
		\label{semi_rot_22}
	\end{subfigure} \\
	\begin{subfigure}[b]{\sfwidth}
		\resizebox{\textwidth}{!}{\input{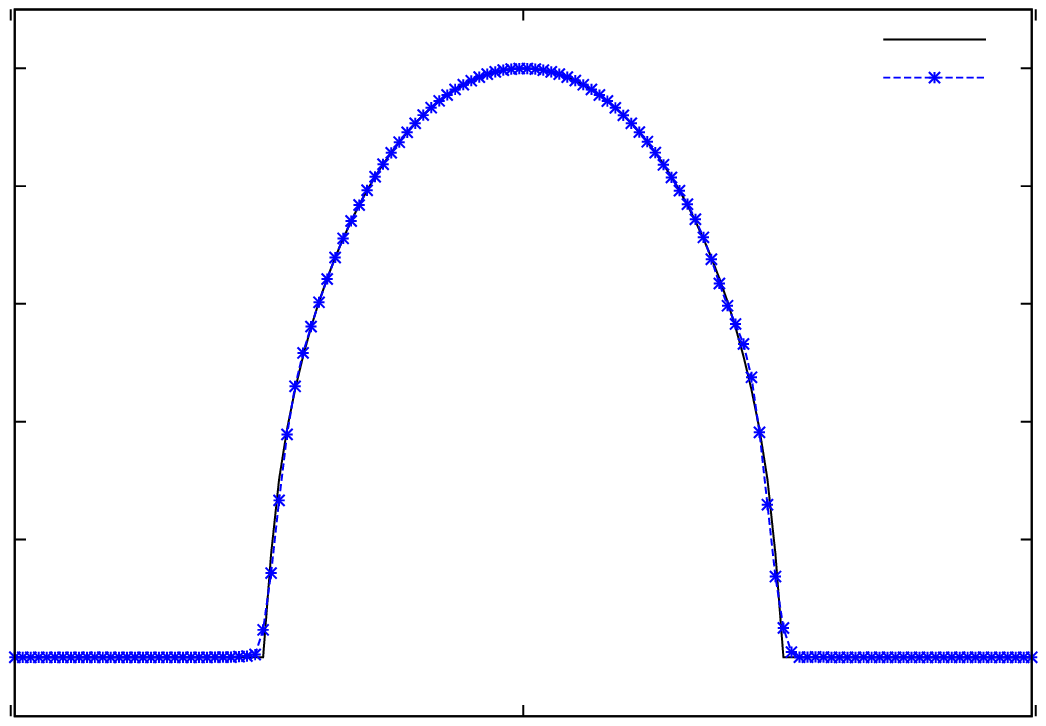}}
		\caption{Centerline Comparison}
		\label{semi_const_23}
	\end{subfigure}
	\begin{subfigure}[b]{\sfwidth}
		\resizebox{\textwidth}{!}{\input{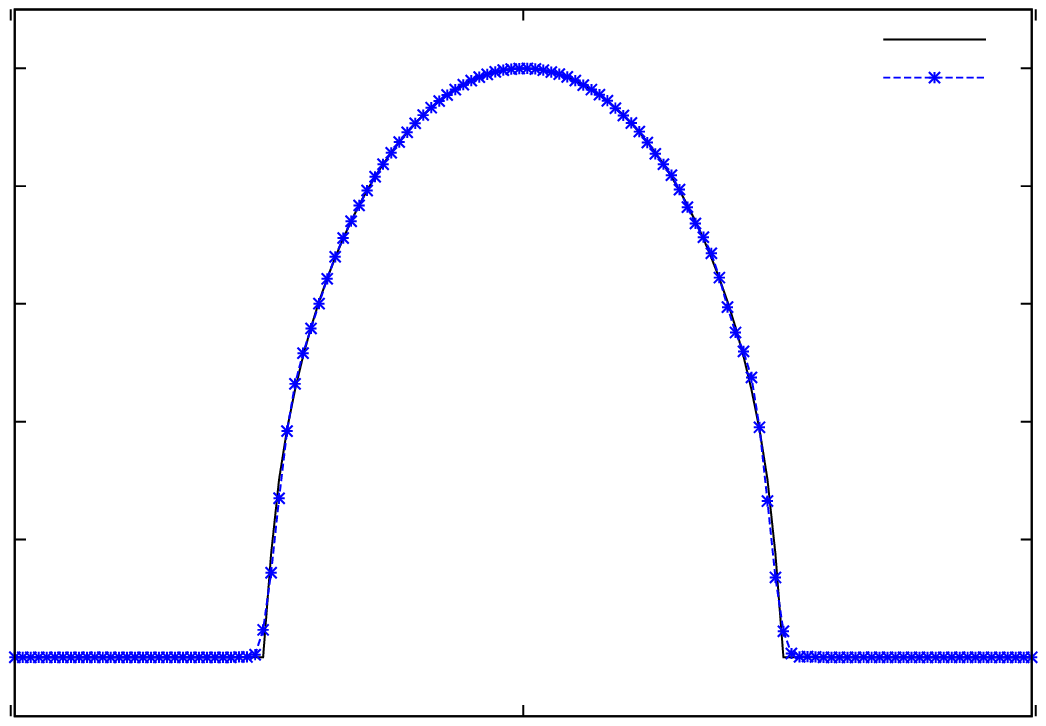}}
		\caption{Centerline Comparison}
		\label{semi_rot_23} 
	\end{subfigure}
\caption[semi_2d]{Semi-ellipse Solutions ($t=1.0,\ N=128$, $\sigma = 0.8$)}
\label{semi_2d}
\end{figure}

\begin{figure} [tb]
\centering
	\begin{subfigure}[b]{\sfwidth}
		\inputImgEps{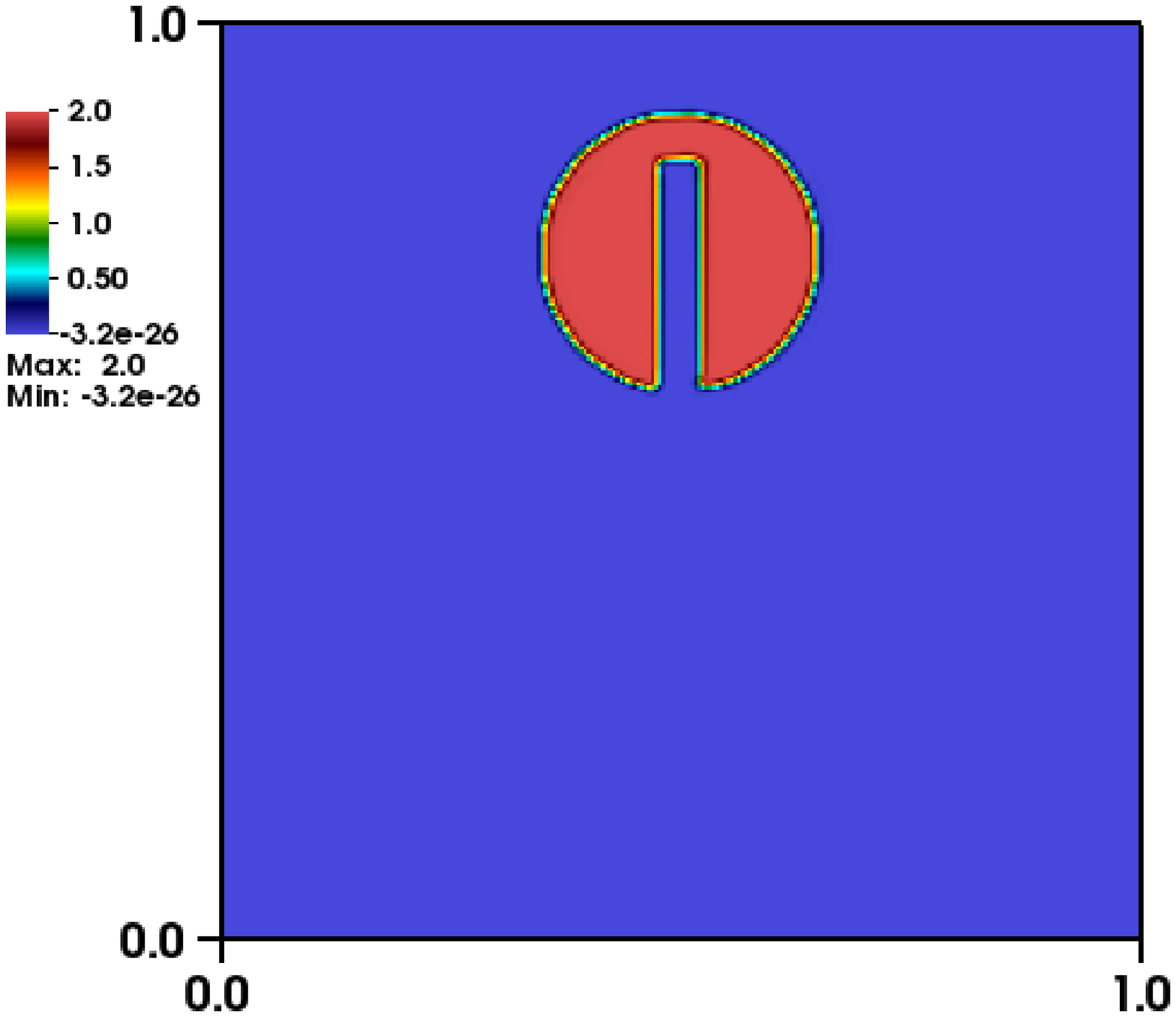}
		\caption{Solid Body Rotation Solution}
		\label{cylinder_22}
	\end{subfigure}
	\begin{subfigure}[b]{\sfwidth}
		\resizebox{\textwidth}{!}{\input{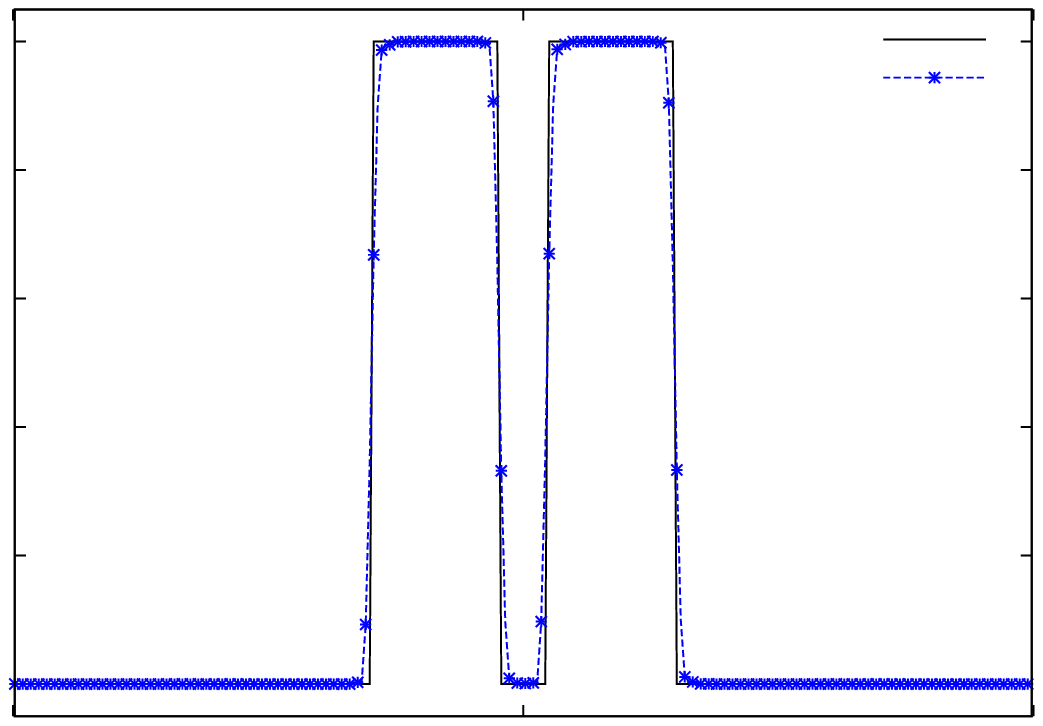}}
		\caption{Centerline Comparison}
		\label{cylinder_23}
	\end{subfigure}
\captionsetup{width=\textwidth}
\caption[cylinder]{Slotted Cylinder under Solid Body Rotation ($t=1.0,\ N=256$, $\sigma = 0.8$)}
\label{cylinder}
\end{figure}

The limiter also performs quite well at representing discontinuities in two dimensions.
Various solution plots for discontinuous initial conditions are presented (\ref{square_2d} - \ref{cylinder}).
All of the two dimensional plots were generated using the ninth order scheme in space, running right at the method's theoretical stability limit ($\sigma = 0.8$ in 2D).
The square solution under constant velocity has few, if any, ripples and is nicely bounded (\ref{square_2d}).
There is some distortion of the corners, particularly at the top-left and bottom-right.
The square solution under solid body rotation looks similar to the constant velocity solution and the corner issue is mitigated.

The semi-ellipse solution is likewise well resolved (\ref{semi_2d}).
As with the one dimensional case there are some dispersive errors on the leading edge, but they are small.
The solutions under both velocity fields are accurately bounded.
The semi-ellipse solution under solid body rotation was centered at $\mb{x}^{\text{solid}}_{c} = \left( 1.0 , 1.5 \right)$
to keep the edge of the conditon away from the domain boundary.
The domain was also expanded to $\mb{x}_{\mb{i}} \in [0,2]$.

The final test was the slotted cylinder (\ref{cylinder}).
The limiter method keeps the solution bounded and resolves the fronts quite nicely.
At lower grid resolutions the slot can fill in and the bounds may not be enforced.
But as the grid is refined, both of these issues are resolved.

\section{Conclusions}

We presented a new flux limiter based upon FCT that retains high-order accuracy for smooth solutions and captures fronts well.
Our algorithm presented here uses CTU for low-order fluxes, upwind schemes for high-order fluxes, and RK4 for time integration.
Our additions to the previous FCT method included a new computation for the extrema and an expanded pre-constraint on the high-order fluxes.
This new pre-constraint is more restrictive than the original one, and seeks to only constrain the fluxes near discontinuities in the solution.

Extensions for this work are applying the limiter to systems of hyperbolic conservation laws, possibly developing new high-order upwind methods with corner-coupling, and further improving the pre-constraint on the high-order fluxes.
Applying the limiter to hyperbolic systems is the most straightforward extension of this work.
It would be relatively simple and informative to apply this limiter in a gas dynamics solver.
High-order, corner-coupled upwind methods for use with general multi-stage time integrators could remove the dimensional dependence of the stability.
However, no upwind methods of this nature currently exists.
The pre-constraint on the high-order fluxes is another area where additional study could pay off.
In this work we found that the pre-constraint affects a delicate balance between effectively representing discontinuities and retaining high-order accuracy in smooth areas. 
It was relatively simple to get one or the other.
Ensuring both required testing many versions of the pre-constraint.
Different versions of the pre-constraint have been proposed \cite{DeVore1998},
but they did not produce good results with the algorithm presented here.

Other versions of FCT \cite{FCT_Book} have implemented high-order centered difference schemes with hyper-diffusive fluxes.
These schemes, like the high-order upwind schemes, seek to add diffusion to the large wave-number solution components.
These centered difference schemes were not examined in this study, but they do offer a possible alternative to the upwind ones used here.

\bibliographystyle{plain}
\bibliography{main}

\end{document}